\documentclass{article}

\usepackage{arxiv}

\usepackage{graphicx}
\usepackage{amsfonts,amsmath,amssymb,mathrsfs,epsfig,bigints,a4,nicefrac,tabularx}
\usepackage{latexsym}
\usepackage{natbib}
\usepackage{enumitem}
\usepackage{setspace}
\usepackage{pstricks}
\usepackage{hyperref,multirow}
\usepackage{lscape}

\usepackage{amsthm}

\numberwithin{equation}{section}
\newcommand{\R}{\mathbb{R}}

\newcommand{\BE}{{\mathbb{E}}}

\newcommand{\RR}{{\mathbb{R}}}

\newcommand{\PP}{{\mathbb{P}}}
\newcommand{\vertk}{\stackrel{\mbox{\scriptsize ${\cal D}$}}{\longrightarrow}}
\newcommand{\edist}{\stackrel{\mbox{\scriptsize ${\cal D}$}}{=}}
\newcommand{\stk}{\stackrel{\mbox{\scriptsize $\PP$}}{\longrightarrow}}
\newcommand{\fsk}{\stackrel{\mbox{\scriptsize a.s.}}{\longrightarrow}}
\newcommand{\LL}{{\textrm L}^2}

\setitemize{leftmargin=5.4mm}

\newtheorem{prop}{Proposition}

\makeatletter
\def\blfootnote{\xdef\@thefnmark{}\@footnotetext}
\makeatother

\title{Tests for multivariate normality -- a critical review with emphasis on weighted $L^2$-statistics}

\author{Bruno Ebner\\
 Institute of Stochastics, \\
Karlsruhe Institute of Technology (KIT), \\
Englerstr. 2, D-76133 Karlsruhe. \\
\texttt{Bruno.Ebner@kit.edu}\\
\And
Norbert Henze\\
Institute of Stochastics, \\
Karlsruhe Institute of Technology (KIT), \\
Englerstr. 2, D-76133 Karlsruhe. \\
\texttt{Norbert.Henze@kit.edu}\\
}

\begin{document}

\date{\today}
\maketitle

%\blfootnote{
%{\em Key words and phrases} Test for multivariate normality; weighted $L^2$-statistic;  affine invariance;  consistency}

\begin{abstract}
This article gives a synopsis on new developments in affine invariant tests for multivariate normality in an i.i.d.-setting,
with special emphasis on asymptotic properties of several classes of weighted $L^2$-statistics. Since weighted $L^2$-statistics
typically have limit normal distributions under fixed alternatives to normality, they open ground for a neighborhood of model validation for normality. The paper also reviews several other invariant tests for this problem, notably the energy test, and it presents the results of a large-scale simulation study. All tests under study are implemented in the accompanying \texttt{R}-package \texttt{mnt}.
\end{abstract}

\section{Introduction}\label{secintro}
Testing for multivariate normality (for short: MVN) is a topic of ongoing interest. A survey of dozens of MVN-tests, including
graphical procedures for assessing multivariate normality, provide \cite{memu04}. The review of \cite{he02} concentrates on
affine invariant and consistent procedures, and the book of \cite{th02} contains a chapter on testing for MVN.

 In a standard setting, let $X,X_1,X_2,\ldots $ be independent identically distributed (i.i.d.) $d$-variate random (column)  vectors, which are defined
on a common probability space $(\Omega,{\cal A},\PP)$. The distribution of $X$ will be denoted by $\PP^X$.
We write  N$_d(\mu,\Sigma)$ for the $d$-variate normal distribution with expectation $\mu$ and
covariance matrix $\Sigma$, and we let
\[
{\cal N}_d := \{{\rm N}_d(\mu,\Sigma): \mu \in \RR^d, \Sigma \ \text{positive definite}\}
\]
denote the class of all non-degenerate $d$-variate normal distributions. Testing for $d$-variate normality means testing the hypothesis
\[
H_0: \PP^X \in {\cal N}_d,
\]
against general alternatives, on the basis of $X_1,\ldots,X_n$. At the outset, it should be stressed that each model
can merely hold approximately in practice. In particular, there can only be approximate normality, in whatever sense.
Consequently, there is the following basic drawback inherent
in any goodness-of-fit test, not only of $H_0$, but also of other families of distributions: If
a level-$\alpha$-test of $H_0$  does not lead to a rejection of $H_0$, the null hypothesis
is by no means ‘validated’ or ‘confirmed’. Presumably, there is merely  not enough evidence to reject it!
A further fundamental point is that there cannot be an optimal test of $H_0$, if one really wants to detect general alternatives.
In this respect, \cite{ja00} shows that the global power function of {\em any} nonparametric test
is flat on balls of alternatives, except for alternatives coming from a finite-dimensional subspace. Thus, loosely speaking,
each test of $H_0$ has its own `non-centrality'.

Regarding the task of reviewing MVN-tests here in 2020, we cite \cite{memu04}, who write `the continuing proliferation of papers
with new methods of assessing MVN makes it virtually impossible for any single survey article to cover all available tests'.
And they continue: `When compared to the amount of work that has been done in developing these tests, relatively little work has been done in
evaluating the quality and power of the procedures'.

This review can also only be partial. We will take the above testing problem seriously and concentrate on genuine tests of $H_0$
that have been proposed since the review \cite{he02}, and we will judge each of these according to the following
points of view:
\begin{itemize}
\item affine invariance
\item theoretical properties (limit distributions under $H_0$ and under fixed and contiguous alternatives to $H_0$, consistency)
\item feasibility with respect to sample size and dimension.
\end{itemize}
Thus, e.g., we will not deal with tests for $H_0$ that allow for $n \le d$ (see \cite{tftw05} or \cite{yahi19}),
since the condition $n \ge d+1$ is necessary to decide whether
the underlying covariance matrix is non-degenerate or not. Moreover, unlike the review of \cite{memu04}, we will not discuss purely graphical procedures, as
proposed in \cite{holg06}. We will also not embark upon a review of tests for normality in non-i.i.d.-settings, like
testing for Gaussianity of the innovations in MGARCH processes (see, e.g.,  \cite{leng11} or \cite{llp14}), or situations with incomplete data (see, e.g., \cite{yhr15}),
since such a task would go beyond the scope of this review. We will also not review tests for Gaussianity in infinite-dimensional Hilbert spaces, see, e.g.,
\cite{ghk20} or \cite{kece19}.

Regarding affine invariance, notice that the class ${\cal N}_d$ is closed with respect to full rank affine transformations.
Hence, any `genuine' statistic $T_n=T_n(X_1,\ldots,X_n)$ (say) for testing $H_0$ should satisfy
$T_n(AX_1+b,\ldots,AX_n+b) = T_n(X_1,\ldots,X_n)$ for each
regular $(d\times d)$-matrix $A$ and each $b \in \RR^d$. Otherwise, it would be possible to reject $H_0$ on given data and
do not object against $H_0$ on the same data, after performing a rotation, which makes little, if any, sense. In the sequel, let
\begin{equation}\label{scaledres}
Y_{n,j} = S_n^{-1/2}(X_j-\overline{X}_n), \quad j=1,\ldots,n,
\end{equation}
denote the so-called {\em scaled residuals}. Here, $\overline{X}_n = n^{-1}\sum_{j=1}^n X_j$ is the sample mean,
$S_n = n^{-1}\sum_{j=1}^n (X_j-\overline{X}_n)(X_j-\overline{X}_n)^\top$ stands for
the sample covariance matrix of $X_1,\ldots,X_n$,  and  the superscript $\top$ denotes transposition of column vectors.
The matrix $S_n^{-1/2}$ is the unique symmetric square root of $S_n^{-1}$. The latter matrix exists almost surely if $n \ge d+1$ and $\PP^X$ is
absolutely continuous with respect to $d$-dimensional Lebesgue measure, see \cite{eape73}. These assumptions will be standing in what follows.
We remark that $S_n$ is sometimes defined with the factor $(n-1)^{-1}$ instead of $n^{-1}$, but this difference does not have
implications for asymptotic considerations. A good account on finite-sample distribution theory of $Y_{n,1}, \ldots, Y_{n,n}$ under $H_0$
is provided by \cite{ta20}.

Affine invariance is achieved if the test statistic $T_n$ is a function of $Y_{n,i}^\top Y_{n,j}, \ i,j \in \{1,\ldots,n\}$, or if
$T_n$ is a function of (only) $Y_{n,1},\ldots,Y_{n,n}$, and $T_n(OY_{n,1}, \ldots,OY_{n,n})= T_n(Y_{n,1},\ldots,Y_{n,n})$
for each orthogonal $(d\times d)$-matrix $O$. If a statistic $T_n$ is affine invariant (henceforth {\em invariant} for the sake of brevity),
the distribution of $T_n$ under the null hypothesis $H_0$ does
not depend on the parameters $\mu$ and $\Sigma$ of the underlying normal distribution. Thus, regarding distribution theory under $H_0$, we can without
loss of generality assume that $\PP^X = \text{N}_d(0,\text{I}_d)$. Here, $0$ is the origin in $\RR^d$, and I$_d$ is the unit matrix of order $d$.
But invariance of a statistic $T_n$ also entails that it is no restriction to assume $\BE X =0$ and $\BE XX^\top = {\rm I}_d$ when studying
the distribution of $T_n$ under an alternative to $H_0$ that satisfies $\BE \|X\|^2 < \infty$, where $\|\cdot\|$ denotes the Euclidean norm in $\RR^d$.

As for the second point, i.e., properties of a test of $H_0$ based on a statistic $T_n$ that go beyond mere simulation results,
there should be a sound rationale for the test, which means that there should be good knowledge of what is estimated by $T_n$ if the underlying distribution
is not normal. This rationale is intimately connected to the property of consistency. If $T_n$ is some invariant statistic, it must be regarded
-- perhaps after some suitable normalization -- as an estimator of some invariant functional ${\cal T}(P)$ of
the unknown underlying distribution $P$, where $P=\PP^X$. This means that ${\cal T}(P)= {\cal T}(\widetilde{P})$ if $\widetilde{P}$ is a full rank
affine image of $P$, whence ${\cal T}(\cdot)$ is constant over the class ${\cal N}_d$. For such a functional, consistency of a test based on ${\cal T}$
against general alternatives can not be expected if ${\cal T}$ does not characterize the class ${\cal N}_d$, in the sense that there are
$P_1 \in {\cal N}_d$ and $P_2 \notin {\cal N}_d$ such that ${\cal T}(P_1) = {\cal T}(P_2)$. Examples of non-characterizing functionals are
time-honored measures of multivariate skewness and kurtosis, see Section \ref{secskurt}. The most prominent of this group of tests is Mardia's invariant non-negative
skewness functional
\begin{equation}\label{mardias}
{\cal T}(P) = \beta_d^{(1)}(P) = \BE\Big{[} \Big{(}(X_1-\mu)^\top \Sigma^{-1}(X_2-\mu)\Big{)}^3\Big{]}.
\end{equation}
Here, $X_1,X_2$ are i.i.d. with distribution $P$, mean $\mu$ and nonsingular covariance matrix $\Sigma$.
The functional $\beta^{(1)}_d$ does not characterize the class ${\cal N}_d$ since it does not only vanish on
${\cal N}_d$, but in particular also for each non-normal elliptically symmetric distribution for which the expectation figuring in \eqref{mardias} exists.
This fact has striking consequences for a standard test of $H_0$ that rejects $H_0$ for large values of the sample counterpart of
$\beta^{(1)}_d$, see Section \ref{secskurt}.

The paper is organized as follows: Section \ref{secmain} gives a thorough account on general aspects of weighted $L^2$-statistics for testing $H_0$, and
besides the class of BHEP-tests, it reviews five recently proposed tests for multivariate normality that are based on either the
characteristic function, the moment generating function, or a combination thereof. Section \ref{sechztest} reviews the Henze--Zirkler test with
bandwidth depending on sample size and dimension, which is not a weighted $L^2$-statistic in the sense of Section \ref{secmain}.
In Section \ref{secenergy}, we summarize the most important features of the meanwhile well established  energy test of \cite{szri05},
 and Section \ref{secpudelko} deals with the test of \cite{pu05}. Section \ref{seccoxsmall} reviews new theoretical results on
 a time-honored test of \cite{cs78}, while Section \ref{secquiroz}  considers the test of \cite{maqu01}, which is based on functions of spherical harmonics.
 In Section \ref{secskurt} we review tests based on skewness and kurtosis, and in Section \ref{secmisc} we try to give a brief account on
 further work on the subject. Section \ref{secsimstud} presents the results of a large scale simulation study that comprises each of the tests
 treated in Sections \ref{secmain} -- \ref{secskurt}. The final Section \ref{secconcl} draws some conclusions, and it gives an outlook for
 further research.

We conclude this section by pointing out some general notation. Throughout the paper, ${\cal B}^d$ stands for the $\sigma$-field of Borel sets
   in $\RR^d$,  ${\mathcal S}^{d-1} := \{x \in \RR^d: \|x\|=1\}$ is the surface of the unit sphere
   in $\RR^d$, and $\Phi(\cdot)$ denotes the distribution function of the standard normal distribution.  The symbol $\vertk$ stands for convergence in distribution of random elements (variables, vectors and processes), and
   $\stk$, $\fsk$  denote convergence in probability and almost sure convergence, respectively. Each limit refers to the setting $n \to \infty$.
The symbol $\edist$ denotes equality in distribution. Throughout the paper, each unspecified integral will be over $\RR^d$. The acronyms (E)MGF and (E)CF
stand for the (empirical) moment generating function and the (empirical) characteristic function, respectively. Finally, we write ${\bf 1}\{A\}$ for the indicator function of an event $A$.

\vspace*{5mm}

\section{Weighted $L^2$-statistics}\label{secmain}
In this chapter, we review the state of the art of weighted $L^2$-statistics for testing $H_0$. These statistics have a long history, and
they are also in widespread use for goodness-of-fit problems with many other distributions, see, e.g., \cite{beh17}. A weighted
$L^2$-statistic for testing $H_0$ takes the form
\begin{equation}\label{defgenerall2}
T_n =  \int Z_n^2(t) w(t)\, {\rm d}t.
\end{equation}
Here, $Z_n(t) = z_n(X_1,\ldots,X_n,t)$, $z_n$ is a real-valued  measurable function defined on
the ($n+1$)-fold cartesian product of $\RR^d$, and $w:\RR^d \to \RR$ is a non-negative weight function satisfying
\[
\int z_n^2(x_1,\ldots,x_n,t) w(t)\, {\rm d}t < \infty  \quad {\textrm{ for each }} (x_1,\ldots,x_n) \in (\RR^d)^n.
\]
The function $z_n$ can also be vector-valued; then $Z_n^2(t)$ in \eqref{defgenerall2} is replaced with $\|Z_n(t)\|^2$.
Typically, $Z_n(t)$   takes the form
\begin{equation}\label{defallgzn}
Z_n(t) = \frac{1}{\sqrt{n}} \sum_{j=1}^n \ell\left(t^\top Y_{n,j}\right), \quad t \in \RR^d,
\end{equation}
where $\ell(\cdot)$ is some measurable function satisfying $\int \BE \big{[} \ell^2(t^\top X)\big{]}w(t) \, {\rm d}t < \infty$, and
$\BE \big{[} \ell(t^\top X)\big{]} =0$, $t \in \RR^d$, if $X \edist {\rm N}_d(0,{\rm I}_d)$. In view of \eqref{defgenerall2}, a natural setting to study asymptotic properties of $T_n$ is the separable Hilbert space
$\mathbb{H} := {\rm L}^2(\RR^d,{\cal B}^d,w(t){\rm d}t)$  of (equivalence classes) of measurabe functions
on $\RR^d$ that are square- integrable with respect to  $w(t)\textrm{d}t$. If $\|f\|_\mathbb{H}:= \left(\int f^2(t) w(t)\, {\rm d}t\right)^{1/2}$
denotes the norm of $f \in \mathbb{H}$, then $T_n = \|Z_n\|^2_\mathbb{H}$. The general approach to derive the limit distribution
of $T_n$ under $H_0$ is to prove $Z_n \vertk Z$ for some centred Gaussian random element of $\mathbb{H}$, whence
$T_n \vertk \|Z\|^2_\mathbb{H}$ by the continuous mapping theorem. To this end, it is indispensable to approximate
$Z_n$ figuring in \eqref{defallgzn} by a suitable random element $Z_{n,0}$ of $\mathbb{H}$ of the form
\begin{equation}\label{defznnu}
Z_{n,0}(t) = \frac{1}{\sqrt{n}} \sum_{j=1}^n \ell_0(t^\top X_j),
\end{equation}
where $\BE[\ell_0(t^\top X)] =0$, $t \in \R^d$, $\int \BE[\ell_0^2(t^\top X)]w(t) {\rm d}t < \infty$, and
$\|Z_n - Z_{n,0}\|_\mathbb{H} \stk 0$. The central limit theorem in Hilbert spaces
(see, e.g., Theorem 2.7 in \cite{bo00}), then yields
$Z_{n,0} \vertk Z$ for some centred Gaussian element of $\mathbb{H}$ having covariance kernel
\[
K(s,t) = \BE \big{[} \ell_0(s,X)\ell_0(t,X)\big{]}, \quad s,t \in \RR^d.
\]

 The distribution of $Z$ is uniquely determined by the kernel $K(\cdot,\cdot)$, and the distribution of $\|Z\|^2_{\mathbb{H}}$
 is that of $\sum_{j=1}^\infty \lambda_j N_j^2$, where the $N_j$ are i.i.d. standard normal random variables, and $\lambda_j$, $j=1,2,\ldots$, are the
 positive eigenvalues corresponding to eigenfunctions $f_j$ of the (linear second-order homogeneous Fredholm) integral equation
\begin{equation}\label{int:eq}
\lambda f(s) = \int K(s,t) f(t) w(t)\,{\rm d}t,\quad s\in \R^d,
\end{equation}
see, e.g., \cite{ks47}. The problem of finding the eigenvalues and associated eigenfunctions of \eqref{int:eq} is  called the {\em kernel eigenproblem}.
 In this respect, hitherto none of the integral equations corresponding to the
 test presented in this section has been solved explicitly. Notice that knowledge of the largest eigenvalue $\lambda_{max}$ (say)
    opens ground for the calculation of the approximate Bahadur slope and hence for statements on the Bahadur efficiency which,
     for asymptotically normal statistics, typically coincides with the Pitman efficiency, for details see \cite{bah60} and \cite{nikitin95}.

To find a random element $Z_{n,0}$ of the form \eqref{defznnu} that approximates $Z_n$, one has to evaluate the effect of
replacing $Y_{n,j}$ in \eqref{defallgzn} with $X_j$. Putting
$\Delta_{n,j} = Y_{n,j}-X_j$, $j=1,\ldots,n$, the following result, taken from \cite{deh19}, is helpful.

\begin{prop}\label{propossumselta}
Let $X, X_1,X_2, \ldots $ be i.i.d. random vectors satisfying $\BE\|X\|^4 < \infty$,  $\BE(X) =0$ and $\BE(X X^\top) = {\rm I}_d$.
We then have
\[
\sum_{j=1}^n \|\Delta_{n,j}\|^2 = O_{\PP}(1), \quad
\frac{1}{n} \sum_{j=1}^n \|\Delta_{n,j}\|^2 \fsk 0, \quad \max_{j=1,\ldots,n} \|\Delta_{n,j}\| = o_\PP\left(n^{-1/4}\right).
\]
\end{prop}
Since $\ell(t^\top Y_{n,j}) = \ell(t^\top X_j + t^\top \Delta_{n,j})$, the function $\ell(\cdot)$ must be smooth enough to allow for a
Taylor expansion.  To tackle the linear part in this expansion, it is crucial to have some information on
$\Delta_{n,j} = (S_n^{-1/2}-{\rm I}_d)X_j - S_n^{-1/2}\overline{X}_n$. Such information is provided by display (2.13) of \cite{hewa97}, according to which
\[
 \sqrt{n}(S_n^{-1/2}- {\rm I}_d) = - \frac{1}{2\sqrt{n}} \sum_{j=1}^n \left(X_jX_j^\top - {\rm I}_d \right) + O_\PP\left(n^{-1/2}\right).
\]
Since Proposition \ref{propossumselta} holds under general assumptions, one may often obtain asymptotic normality of
weighted $L^2$-statistics under fixed alternatives. To this end, notice  that
\[
\frac{T_n}{n} =   \int \bigg{(} \frac{1}{n}\sum_{j=1}^n \ell(t^\top Y_{n,j})\bigg{)}^2  w(t)\, {\rm d}t.
\]
Under suitable conditions, we will have $T_n/n \stk \Delta$, where $\Delta = \|z\|^2_\mathbb{H}$, and $z(t) = \BE \big{[} \ell(t^\top X)\big{]}$,
$z \in \RR^d$. An immediate consequence of this stochastic convergence is the consistency of a test for $H_0$ based on $T_n$ against
each alternative distribution that satisfies $\Delta >0$. But we have more! Writing $\langle u,v \rangle_\mathbb{H} = \int u(t)v(t) w(t) \, \text{d}t$ for the inner product in
$\mathbb{H}$, there is the decomposition
\begin{eqnarray*}
\sqrt{n}\left(\frac{T_n}{n} - \Delta\right)  & = & \sqrt{n} \left(\|Z_n\|^2_{\mathbb{H}} - \|z\|^2_{\mathbb{H}} \right)\\
& = &  \sqrt{n}\big{\langle} Z_n-z,Z_n+z \big{\rangle}_\mathbb{H}\\
& = &  \sqrt{n}\big{\langle} Z_n-z,2z + Z_n-z \big{\rangle}_\mathbb{H}\\
& = &  2\big{\langle} \sqrt{n}(Z_n-z),z\big{\rangle}_\mathbb{H} + \frac{1}{\sqrt{n}}\|\sqrt{n}(Z_n-z)\|^2_{\mathbb{H}}.
\end{eqnarray*}
These lines carve out the quintessence of asymptotic normality of weighted $L^2$-statistics under fixed alternatives. Namely, if one can show that the sequence
$V_n := \sqrt{n}(Z_n-z)$ of random elements of $\mathbb{H}$ converges in distribution to some centred Gaussian random element $V$ of
$\mathbb{H}$, then, by the continuous mapping theorem and Slutski's lemma, we have
\begin{equation}\label{asynfix}
\sqrt{n}\left(\frac{T_n}{n} - \Delta\right) \ \vertk \ {\rm{N}}(0,\sigma^2),
\end{equation}
where
\[%\begin{equation}\label{sigmaqu}
\sigma^2 = 4 \iint K(s,t) z(s)z(t) w(s)w(t) \ \text{d}s\text{d}t,
\]%\end{equation}
and $K(\cdot,\cdot)$ is the covariance kernel of $V$, see Theorem 1 of \cite{beh17}. As a consequence, if $\widehat{\sigma}_n^2$ is a consistent estimator of
$\sigma^2$ based on $X_1,\ldots,X_n$, then, for given $\alpha \in (0,1)$,
\begin{equation}\label{confintbhep}
I_{n,1-\alpha} = \bigg{[} \frac{T_n}{n} - \Phi^{-1}\left(1\! - \! \frac{\alpha}{2}\right) \frac{\widehat{\sigma}_{n}}{\sqrt{n}},
                               \frac{T_n}{n} + \Phi^{-1}\left(1\! - \! \frac{\alpha}{2}\right) \frac{\widehat{\sigma}_{n}}{\sqrt{n}}\bigg{]}
\end{equation}
is an asymptotic confidence interval for $\Delta$ of level $1-\alpha$. Moreover, from \eqref{asynfix} and Slutski's lemma, we have
\begin{equation}\label{asyss}
\frac{\sqrt{n}}{\widehat{\sigma}_n}\left(\frac{T_n}{n} - \Delta\right) \ \vertk \ {\rm{N}}(0,1),
\end{equation}
which opens the ground for a validation of a certain neighborhood of $H_0$. Namely, suppose that we want to tolerate a given `distance' $\Delta_0$
to the class ${\cal N}_d$. We may then consider the `inverse' testing problem
\[
H_{\Delta_0}: \Delta(\PP^X) \ge \Delta_0 \  \textrm{ against } \ K_{\Delta_0}: \Delta(\PP^X) < \Delta_0.
\]
Here, the dependence of $\Delta$ on the underlying distribution $\PP^X$ has been made explicit.

From \eqref{asyss}, the test which rejects $H_{\Delta_0}$ if
\[
\frac{T_n}{n} \ \le \ \Delta_0 - \frac{\widehat{\sigma}_n}{\sqrt{n}} \Phi^{-1}(1-\alpha),
\]
has asymptotic level $\alpha$, and it is consistent against general alternatives, see
Section 3.3 of \cite{beh17}. Notice that this test is in the spirit of bioequivalence testing
 (see, e.g., \cite{cfm07}, \cite{demu}  or \cite{wellek}), since it aims at validating a certain neighborhood of a hypothesized model.

We now review the time-honored class of BHEP-tests and several recently suggested $L^2$-statistics for testing $H_0$. Each of these statistics
has an upper rejection region, and it is invariant, because it is a function of $Y_{n,j}^\top Y_{n,k}$, where $j,k \in \{1,\ldots,n\}$.

\subsection{The BHEP-tests}\label{secbhep}
Generalizing a test for
univariate normality based on the ECF due to \cite{ep83},
the first proposals for weighted  $L^2$-statistics for testing $H_0$ are due to \cite{bahe88}  and \cite{hezi90}, who considered the
statistic
\begin{equation}\label{defBHEP}
\text{BHEP}_{n,\beta} = n \int \big{|} \Psi_n(t) - \Psi_0(t)\big{|}^2 w_\beta(t) \, \text{d}t.
\end{equation}
Here,
\begin{equation}\label{defpsint}
\Psi_n(t) = \frac{1}{n} \sum_{j=1}^n \exp(\text{i}t^\top Y_{n,j}), \qquad t \in \RR^d,
\end{equation}
denotes the ECF of $Y_{n,1}, \ldots, Y_{n,n}$, $\Psi_0(t) = \exp(-\|t\|^2/2)$ is the CF of
the distribution ${\rm N}_d(0,{\rm I}_d)$,
and the weight function $w_\beta$ is given by
\begin{equation}\label{weightfunction}
w_\beta(t) = \left(2\pi \beta^2\right)^{-d/2} \exp\left(- \frac{\|t\|^2}{2\beta^2}\right),
\end{equation}
where $\beta >0$ is a fixed constant. That $\text{BHEP}_{n,\beta}$ is indeed of the type \eqref{defgenerall2}
will become clear from the representation  \eqref{dstbhepl2}.

Whereas \cite{bahe88} studied the special case $\beta =1$, the general case
was treated by \cite{hezi90}.
An extremely
appealing feature of the weight function $w_\beta$ in \eqref{weightfunction} is that BHEP$_{n,\beta}$ takes the feasible form
\begin{eqnarray}\label{BHEP}
\text{BHEP}_{n,\beta} & = & \frac{1}{n} \sum_{j,k=1}^n \exp\left(-\frac{\beta^2\|Y_{n,j}-Y_{n,k}\|^2}{2}\right)\\ \nonumber
& & \quad - \frac{2}{(1+\beta^2)^{d/2}} \sum_{j=1}^n \exp\left(- \frac{\beta^2\|Y_{n,j}\|^2}{2(1+\beta^2)} \right) + \frac{n}{(1+2\beta^2)^{d/2}}.
\end{eqnarray}
The BHEP-test is the most thoroughly studied class of tests for multivariate normality.
\cite{cs89} coined the acronym BHEP for this class of tests for $H_0$, after early developers of the idea, and
he proved that $\liminf_{n\to \infty} n^{-1} {\rm BHEP}_{n,\beta} \ge C(\PP^X,\beta) >0$ almost surely  for some constant $C(\PP^X,\beta)$ if
$\PP^X$ does not belong to ${\cal N}_d$. As a consequence, a test for normality based on ${\rm BHEP}_{n,\beta}$ is
consistent against any alternative.

If $\BE \|X\|^2 < \infty$ and $\BE X  =0$, $\BE XX^\top  = {\rm I}_d$ (the last two assumptions entail no
loss of generality in view of invariance), then
\begin{equation}\label{fskbhep}
\frac{1}{n} {\rm BHEP}_{n,\beta} \fsk \Delta_\beta := \int \big{|} \Psi(t) - \Psi_0(t) \big{|}^2 w_\beta(t)\, {\rm d}t
\end{equation}
(\cite{bahe88}), where $\Psi(t) = \BE \exp({\rm i}t^\top X)$, $t \in \R^d$, is the CF of $X$. Hence, $\Delta_\beta = \Delta_\beta(\PP^X)$ is the
functional associated with the BHEP-test. Using a Hilbert space setting,  \cite{gu00} proved \eqref{asynfix}
for $T_n = {\rm BHEP}_{n,\beta}$, where $\Delta = \Delta_\beta$ and $\sigma^2 = \sigma_\beta^2$ depend on $\beta$,
under each alternative distribution satisfying $\BE\|X\|^4 < \infty$.
Moreover, \cite{gu00} obtained a sequence $\widehat{\sigma}^2_{n,\beta}$ of consistent estimators of $\sigma_\beta^2$ and thus
an asymptotic confidence interval of the type \eqref{confintbhep}.

In view of the representation \eqref{BHEP},
\cite{bahe88} and \cite{hezi90} obtained the limit null distribution of BHEP$_{n,\beta}$ as $n \to \infty$ my means of  the theory of V-statistics with estimated parameters.
Upon observing that
\begin{equation}\label{dstbhepl2}
\text{BHEP}_{n,\beta} = \int Z_n^2(t) \, w_\beta(t)\, \text{d}t,
\end{equation}
where $Z_n(t) = n^{-1/2} \sum_{j=1}^n \left( \cos(t^\top Y_{n,j}) + \sin(t^\top Y_{n,j}) - \Psi_0(t)\right)$,
\cite{hewa97} considered $Z_n(\cdot)$ as a random element in a certain Fr\'{e}chet space of random functions, and they showed that
$Z_n$ converges in distribution in that space to some centred Gaussian random element $Z$, see Theorem 2.1 of \cite{hewa97}. Moreover, BHEP$_{n,\beta} \vertk \int Z^2(t)\, w_\beta(t) \, \text{d} t$, and the
test is able to detect a sequence of contiguous alternatives that approach $H_d$ at the rate $n^{-1/2}$.
\cite{hewa97} also obtained the first three moments of the limit null distribution of BHEP$_{n,\beta}$. Finally, the class of BHEP-tests is `closed at the boundaries' $\beta \to 0$ and $\beta \to \infty$ since,
elementwise on the underlying probability space, we have
\begin{equation}\label{bhepbeta0}
\lim_{\beta \to 0} \frac{{\rm BHEP}_{n,\beta}}{\beta^6} = \frac{n}{6} \cdot b_{n,d}^{(1)} + \frac{n}{4} \cdot \widetilde{b}_{n,d}^{(1)},
\end{equation}
where $b_{n,d}^{(1)}$ and $\widetilde{b}_{n,d}^{(1)}$ are given in \eqref{defskewnkurtma} and \eqref{defskmorosz}, respectively, see \cite{he97b}.
Thus, as $\beta \to 0$, a scaled version of ${\rm  BHEP}_{n,\beta}$ is approximately a linear combination of two measures of multivariate skewness.
The limit distribution of the right-hand side of \eqref{bhepbeta0} under general distributional assumptions on $X$ has been studied by \cite{he97b}.
 Last but not least,
we have
\begin{equation}\label{bhepbetainf}
\lim_{\beta \to \infty} \beta^d \left({\rm BHEP}_{n,\beta}-1 \right) = \frac{n}{2^{d/2}} - 2 \sum_{j=1}^n \exp\left( - \frac{\|Y_{n,j}\|^2}{2}\right),
\end{equation}
see \cite{he97b}.
Hence, as $\beta \to \infty$, rejection of $H_0$ for large values of ${\rm BHEP}_{n,\beta}$ means rejection of $H_0$ for {\em small} values of $\sum_{j=1}^n \exp(-\|Y_{n,j}\|^2/2)$.
The latter statistic, like Mardia's measure of multivariate kurtosis $b_{n,d}^{(2)}$ (see \eqref{defskewnkurtma}), merely investigates an aspect of the `radial part' of the underyling distribution.

Guided by theoretical and simulation based results in the univariate case, \cite{tenr09} performed an extensive simulation
study on the power of the BHEP test for dimensions $d \in \{2,3,\ldots,10,12,15 \}$ and sample sizes $n \in \{20,40,60,80,100\}$.
He concluded that the choice $\beta_n$ given in \eqref{betahz} gives `the best results for long tailed or moderately skewed
alternatives, but it also produces very poor results for short tailed alternatives'. If no relevant information about
the tail of the alternatives is available, he strongly recommends the use of
$\beta = \sqrt{2}/(1.376+0.075d)$
(in fact,his recommendation is in terms of $h = 1/(\beta \sqrt{2})$)), and there are similar recommendations
for short tailed alternatives and long tailed or
moderately skewed alternatives, respectively.

\subsection{A weighted $L^2$-statistic via the moment generating function}\label{secmgfu}
\cite{hejg19} generalized results of  \cite{heko20} to the multivariate case and considered a MGF analogue
to the BHEP-test statistic. Letting
 \begin{equation}\label{EMGF}
 M_n(t) = \frac{1}{n} \sum_{j=1}^n \exp\left(t^\top Y_{n,j}\right), \quad t\in \RR^d,
 \end{equation}
denote the EMGF of $Y_{n,1}, \ldots,Y_{n,n}$, and writing
 $ M_0(t)=\exp(\|t\|^2/2)$, $t \in \RR^d$,  for the  MGF of the standard normal distribution N$_d(0,\textrm{I}_d)$, the test statistic
is
\begin{equation}\label{teststatHJM}
{\rm HJ}_{n,\gamma} = n \int \left(M_n(t) - M_0(t) \right)^2 \, \widetilde{w}_\gamma(t) \, \textrm{d}t,
\end{equation}
where
\begin{equation}\label{weightbeta}
\widetilde{w}_\gamma(t) = \exp\left(-\gamma \|t\|^2 \right),
\end{equation}
and $\gamma >2$ is some fixed parameter. Notice that the condition $\gamma>1$  is necessary for the integral in (\ref{teststatHJM}) to be finite, and
the more stringent condition $\gamma >2$ is needed for asymptotics under $H_0$. The test statistic ${\rm HJ}_{n,\gamma}$ has a representation analogous to \eqref{BHEP} (see display (1.4) of \cite{hejg19}). Elementwise on the underlying probability space, we have
\begin{equation}\label{limitgainf}
\lim_{\gamma \to \infty} \gamma^{3+d/2} \, \frac{6 {\rm HJ }_{n,\gamma}}{\pi^{d/2}} = \frac{n}{6} \cdot b_{n,d}^{(1)} +  \frac{n}{4} \widetilde{b}_{n,d}^{(1)}
\end{equation}
which, interestingly, is the same limit as in \eqref{bhepbeta0}.
By working in the Hilbert space ${\rm L}^2(\RR^d,{\cal B}^d,\widetilde{w}_\gamma(t){\rm d}t)$  of (equivalence classes) of measurabe functions
on $\RR^d$ that are square- integrable with respect to  $\widetilde{w}_\gamma(t)\textrm{d}t$,
\cite{hejg19} derived the limit null distribution of  ${\rm HJ}_{n,\gamma}$, which is that of
${\rm HJ}_{\infty,\gamma} := \int W^2(t) \widetilde{w}_\gamma(t) \, {\rm d}t$, where
$W$ is some centred Gaussian random element of that space.
\cite{hejg19} also obtained the expectation and the variance of ${\rm HJ}_{\infty,\gamma}$. Moreover, if
$X$ is a (standardized)  alternative distribution with the property $M(t) :=$ $\BE(\exp(t^\top X)) < \infty$, $t \in \RR^d$, then
\begin{equation}\label{liminfhj}
\liminf_{n \to \infty} \frac{{\rm HJ}_{n,\gamma}}{n} \ \ge \  \int \left(M(t) -M_0(t)  \right)^2 \, \widetilde{w}_\gamma(t) \, \rm{d} t \qquad \PP\text{-almost surely}.
\end{equation}
This inequality implies the consistency of the MVN test based on ${\rm HJ}_{n,\gamma}$ against those alternatives that have a finite MGF.
Indeed, one may conjecture that this test is consistent against {\em any} alternative to $H_0$.

\subsection{A test based on a characterization involving the MGF and the CF}\label{sececfmgf}
\cite{vo14} proved a characterization of the univariate  centred normal distribution, which involves both the CF and the MGF.
\cite{hjmm19} generalized this result as follows: If $X$ is a centred $d$-variate non-degenerate
random vector with MGF
 $M(t) = \BE[\exp(t^\top X)]<\infty$,  $t \in \RR^d$, and
$R(t) := \BE[\cos(t^\top X)]$ denotes the real part of the CF of $X$, then
\begin{equation} \label{RCF}
R(t) \, M(t)-1 \ =\ 0 \quad  \textrm{for each } t  \in \RR^d
\end{equation}
holds true if  and only if  $X$ follows some zero-mean normal distribution.

Since $Y_{n,1}, \ldots, Y_{n,n}$ provide an empirical
standardization of $X_1,\ldots,X_n$, a natural test statistic based on \eqref{RCF} is
\[
{\rm HJM}_{n,\gamma} \ := \ n \int \left(R_n(t)\, M_n(t) - 1 \right)^2\, \widetilde{w}_\gamma(t) \, \textrm{d} t,
\]
where
\[%\begin{equation}\label{ECT}
R_n(t) \ := \ \frac{1}{n} \sum_{j=1}^n \cos\left(t^\top Y_{n,j}\right), \quad t \in \RR^d,
\]%\end{equation}
is the  empirical cosine transform of the scaled residuals,  and $M_n(t)$ and $\widetilde{w}_\gamma(t)$ are
given in \eqref{EMGF} and \eqref{weightbeta}, respectively. There is a representation of
${\rm HJM}_{n,\gamma}$ similar to  \eqref{BHEP}, but involving a fourfold sum (see display (3.7) of \cite{hjmm19}).
The main results about ${\rm HJM}_{n,\gamma}$ are as follows: Elementwise on the underlying probability space, we have
\[
\lim_{\gamma \to \infty} \gamma^{3+d/2} \, \frac{8 {\rm HJM}_{n,\gamma}}{\pi^{d/2}} = \frac{n}{6} \cdot b_{n,d}^{(1)} +  \frac{n}{4} \cdot \widetilde{b}_{n,d}^{(1)}.
\]
Interestingly, this is the same linear combination of two measures of skewness as in \eqref{bhepbeta0} and \eqref{limitgainf}.
If $\gamma >1$, then the limit null distribution of ${\rm HJM}_{n,\gamma}$ is that of ${\rm HJM}_{\infty,\gamma} := \int W^2(t) \widetilde{w}_\gamma(t) \, {\rm d}t$,
where $W$ is a centred random element of the Hilbert space ${\rm L}^2(\RR^d,{\cal B}^d,\widetilde{w}(t){\rm d}t)$
with a covariance kernel given in Theorem 5.1 of \cite{hjmm19}. Moreover, that paper also states a formula for
$\BE[{\rm HJM}_{\infty,\gamma}]$ and obtains the inequality
\begin{equation}\label{consistin}
\liminf_{n \to \infty} \frac{{\rm HJM}_{n,\gamma}}{n} \ \ge \  \int \left( R(t)M(t) -1 \right)^2 \, w_\gamma(t) \, \rm{d} t \quad \PP\text{-almost surely},
\end{equation}
which is analogous to \eqref{liminfhj}. We conjecture that also the MVN test based on ${\rm HJM}_{n,\gamma}$ is consistent against
{\em any} non-normal alternative distribution.

\subsection{A test based on a system of partial differential equations for the MGF}%The test of Henze and Visagie}
The novel idea of \cite{hevi19} for constructing a test of $H_0$ is the following:
Suppose that the MGF $M(t) = \mathbb{E} [\exp(t^\top X)]$ of a random vector $X$ exists for each
$t \in \mathbb{R}^d$ and satisfies the system of partial differential  equations
\begin{equation}\label{partial}
\frac{\partial M(t)}{\partial t_j}  = t_j M(t),  \quad t=(t_1,\ldots,t_d)^\top \in \mathbb{R}^d, \quad j=1,\ldots,d.
\end{equation}
Since $M(0) =1$, it is easily seen that the only solution to \eqref{partial} is  $M_0(t) = \exp(\|t\|^2/2)$, $t \in \mathbb{R}^d$,
 which is the MGF of N$_d(0,\textrm{I}_d)$. If $H_0$ holds, the scaled residuals $Y_{n,1},\ldots,Y_{n,n}$ should be approximately independent,
with a distribution close to N$_d(0,\textrm{I}_d)$, at least for large $n$.
Hence, a natural approach for testing $H_0$ is to consider the EMGF $M_n$
of $Y_{n,1},\ldots,Y_{n,n}$, defined in \eqref{EMGF}, and to employ the
weighted $L^2$-statistic
\[%\begin{equation}\label{teststat}
{\rm HV}_{n,\gamma} := n \int \|\nabla M_n(t) - t M_n(t)\|^2 \, \widetilde{w}_\gamma(t)  \, \textrm{d}t,
\]%\end{equation}
where $\nabla f$ stands for the gradient of a function $f:\RR^d \to \RR$, and $\widetilde{w}_\gamma$ is given in \eqref{weightbeta}.
Putting $Y_{n,j,k}^+ = Y_{n,j} + Y_{n,k}$, ${\rm HV}_{n,\gamma}$ takes the feasible form
\[
{\rm HV}_{n,\gamma} \! = \! \frac{1}{n} \! \! \left(\frac{\pi}{\gamma}\right)^{d/2} \! \! \sum_{j,k=1}^n  \! \exp \! \left(\! \frac{\|Y_{n,j,k}^+\|^2}{4 \gamma} \! \right) \! \! \! \left(\! Y_{n,j}^\top Y_{n,k} \! - \!  \frac{\|Y_{n,j,k}^+\|^2}{2\gamma}
\! + \!  \frac{d}{2\gamma} \! + \!  \frac{\|Y_{n,j,k}^+\|^2}{4\gamma^2} \! \right) \! \!.
\]
To derive the limit null distribution of ${\rm HV}_{n,\gamma}$, put
$W_n(t) := \sqrt{n} \left(\nabla M_n(t) - t M_n(t)\right)$. Since
$W_n(t)$ is $\RR^d$-valued, \cite{hevi19}  consider the Hilbert space $\mathbb{H}$, which is the $d$-fold (orthogonal) direct sum $\mathbb{H} := \LL \oplus \cdots \oplus \LL$,
where $\LL = {\rm L}^2(\RR^d,{\cal B}^d,\widetilde{w}(t){\rm d}t)$.
If $\gamma >2$, there is some centred Gaussian random element $W$ of $\mathbb{H}$ with a
covariance {\rm{(}}matrix{\rm{)}} kernel given in display (11) of \cite{hevi19},
so that $W_n \vertk W$ as $n \to \infty$. By the continuous mapping theorem, we then have ${\rm HV}_{n,\gamma} \vertk
  {\rm HV}_{\infty,\gamma} := \int \|W(t)\|^2 \, \widetilde{w}_\gamma(t)   \, {\rm{d}}t$.
\cite{hevi19} also obtain a closed form expression for
$\BE[T_{\infty,\gamma}]$. Moreover, if the MGF $M(t)$ of $X$ exists for each $t \in \RR^d$ and $X$ is standardized, we have
\[%\begin{equation}\label{righthside}
\liminf_{n\to \infty} \frac{{{\rm HV}}_{n,\gamma}}{n} \ge \int \|\nabla M'(t) - tM(t)\|^2 \, \widetilde{w}_\gamma(t) \, {\rm{d}}t \quad \PP\text{-almost\ surely},
\]%\end{equation}
which parallels \eqref{liminfhj} and \eqref{consistin}.

\subsection{A test based on the harmonic oscillator in characteristic function spaces}
\cite{deh19} noticed that the CF $\Psi_0(t) = \exp(-\|t\|^2/2)$ of the distribution
 ${\rm N}(0,{\rm I}_d)$ is the unique solution of the partial differential equation
 \begin{equation}\label{PDEHO}
 \Delta f(x) - (\|x\|^2-d) f(x)=0
 \end{equation} subject to $f(0) =1$, where $\Delta$ is the Laplace operator, see Theorem 1 of \cite{deh19}.
 The operator $-\Delta + \|x\|^2 - d$ is called the {\em harmonic oscillator}, which is a special case of a Schr\"{o}dinger operator.
A suitable statistic for testing $H_0$ that reflects this characterization is
\begin{eqnarray}\label{defdeh}
{\rm DEH}_{n,\gamma} &=&n\int_{\R^d}\left|\Delta\Psi_n(t)- \Delta \Psi_0(t)\right|^2 \widetilde{w}_\gamma(t)\mbox{d}t\\ \nonumber
&=& n\int\bigg{|} \frac{1}{n} \sum_{j=1}^n \|Y_{n,j}\|^2 \exp({\rm i}t^\top Y_{n,j})+(\|t\|^2-d)\Psi_0(t)\bigg{|}^2\widetilde{w}_\gamma(t)\, \mbox{d}t,
\end{eqnarray}
where $\widetilde{w}_\gamma$ is given in \eqref{weightbeta} and $\gamma >0$. The test statistic has the feasible form
\begin{eqnarray*}
{\rm DEH}_{n,\gamma}&=&\left(\frac\pi{\gamma}\right)^\frac{d}{2}\frac1n\sum_{j,k=1}^n\|Y_{n,j}\|^2\|Y_{n,k}\|^2\exp\left(-\frac1{4\gamma}\|Y_{n,j}-Y_{n,k}\|^2\right)\\
&&-\frac{2(2\pi)^{\frac{d}2}}{(2\gamma+1)^{2+\frac{d}{2}}}\sum_{j=1}^n\|Y_{n,j}\|^2\left(\|Y_{n,j}\|^2+2d\gamma(2\gamma+1)\right)\exp\left(-\frac12\frac{\|Y_{n,j}\|^2}{2\gamma+1}\right)\\
&&+n\frac{\pi^{\frac{d}2}}{(\gamma +1)^{2+\frac{d}{2}}}\left(\gamma(\gamma +1)d^2+\frac{d(d+2)}{4}\right).
\end{eqnarray*}
Like the class of BHEP-tests, also the class of tests based on ${\rm DEH}_{n,\gamma}$ is closed at the boundaries $\gamma \to 0$ and $\gamma \to \infty$, since -- elementwise
on the underlying probability space -- we have
\[
\lim_{\gamma \to 0} \left(\frac{\gamma}{\pi}\right)^{d/2} {\rm DEH}_{n,\gamma} = b_{n,d}^{(2)}, \quad
\lim_{\gamma \rightarrow\infty} \frac2{n\pi^{\frac d2}}\gamma^{\frac d2+1} {\rm DEH}_{n,\gamma} =  \widetilde{b}_{n,d}^{(1)}.
\]
Here, $b_{n,d}^{(2)}$ is multivariate kurtosis in the sense of \cite{ma70}, defined in \eqref{defskewnkurtma}, and
$\widetilde{b}_{n,d}^{(1)}$ is skewness in the sense of \cite{morosz93}, see \eqref{defskmorosz}.
\cite{deh19}  proved a Hilbert space central limit theorem for the sequence of random elements
\[
V_n(t) = \frac{1}{\sqrt{n}} \sum_{j=1}^n \left(\|Y_{n ,j}\|^2\big{\{} \cos(t^\top Y_{n,j}) + \sin(t^\top Y_{n,j})\big{\}} - \mu(t)\right), \quad t \in \RR^d,
\]
where $\mu(t) = \BE[\|X\|^2(\cos(t^\top X + \sin(t^\top X))]$, and $X$ is a standardized random vector satisfying $\BE \|X\|^4 < \infty$.
Since $\mu(t) = -\Delta \Psi_0(t)$ if $X \edist {\rm N}_d(0,{\rm I}_d)$ and  ${\rm DEH}_{n,\gamma} = \int V_n^2(t) \widetilde{w}_\gamma(t) \, {\rm d}t$
for that choice of $\mu(t)$, the authors obtained the limit distribution of ${\rm DEH}_{n,\gamma}$ under $H_0$ as well as under
contiguous and fixed alternatives to $H_0$. Under $H_0$, we have ${\rm DEH}_{n,\gamma} \vertk \int V^2(t) \widetilde{w}_\gamma(t) \, {\rm d}t$, where
$V$ is the centred limit Gaussian random element of the sequence $(V_n)$ (with $\mu(t) = - \Delta \Psi_0(t)$). Under contiguous
alternatives that approach $H_0$ at the rate $n^{-1/2}$, the limit distribution of   ${\rm DEH}_{n,\gamma}$ is that of
$\int (V(t)+c(t))^2 \widetilde{W}_\gamma(t) \, {\rm d}t$, where $c(\cdot)$ is a shift function (see Section 6 of \cite{deh19}).
Under a fixed (and because of invariance without loss of generality standardized) alternative distribution satisfying $\BE\|X\|^4  < \infty$, we have
\[
\frac{{\rm DEH}_{n,\gamma}}{n} \rightarrow {\rm D}_\gamma := \int \big{|}\Delta \Psi(t) - \Delta \Psi_0(t)\big{|}^2 \widetilde{w}_\gamma(t) \, {\rm d}t \quad \text{$\PP$-almost surely,}
\]
where $\Psi$ is the CF of $X$. Moreover, the limit distribution of $\sqrt{n}({\rm DEH}_{n,\gamma}/n - {\rm D}_\gamma)$ is a centred normal distribution with a
variance that, under the stronger condition $\BE \|X\|^6 < \infty$,  can be consistently estimated from the data. Thus, by analogy with \eqref{confintbhep}, an asymptotic confidence interval
for ${\rm D}_\gamma$ is available. Notice that, when compared with \eqref{fskbhep}, the almost sure limits above are `Laplacian analogues' of \eqref{fskbhep}.

\subsection{A test  based on a double estimation in a characterizing PDE}
\cite{deh19b} suggested to replace {\em both} of the functions $f$ occurring in   \eqref{PDEHO} by the
ECF $\Psi_n$. Since, under $H_0$, $\Delta \Psi_n(t)$ and $(\|t\|^2 -d)\Psi_n(t)$ should be close to each other for large $n$,
it is tempting to see what happens if, instead of ${\rm DEH}_{n,\gamma}$ defined in   \eqref{defdeh},
we base a test of $H_0$ on the weighted $L^2$-statistic
\[%\begin{equation}\label{defuna}
{\rm DEH}^*_{n,\gamma} = n\int \left|\Delta\Psi_n(t)-\left(\|t\|^2-d\right)\Psi_n(t)\right|^2\widetilde{w}_\gamma(t)\, \mbox{d}t
\]%\end{equation}
and reject $H_0$ for large values of ${\rm DEH}^*_{n,\gamma}$. Putting $D^2_{n,j,k} := \|Y_{n,j}-Y_{n,k}\|^2$,
$E_{n,j,k} = \exp(-D^2_{n,j,k}/(4\gamma))$, $a_{d,\gamma} =2\gamma d(2\gamma\! -\! 1)$,
$b_{d,\gamma} = 16d^2\gamma^3(\gamma\! -\! 1) + 4d(d\! +\! 2)\gamma^2$, $c_{d,\gamma} = (\pi/\gamma)^{d/2}$, and $e_{d,\gamma} = 8d\gamma^2\! -\!  4(d\! +\! 2)\gamma$,
the statistic ${\rm DEH}^*_{n,\gamma}$ has the feasible representation
\begin{eqnarray*}
    {\rm DEH}^*_{n,\gamma}  & = & \frac{c_{d,\gamma}}{n} \!
    \sum_{j,k=1}^n \! \Biggl[ \! \|Y_{n,j}\|^2\|Y_{n,k}\|^2 E_{n,j,k} \! - \! \frac{\|Y_{n,j}\|^2\! + \!  \|Y_{n,k}\|^2}{4\gamma^2}\bigl(D^2_{n,j,k}\! +\!  a_{d,\gamma}\bigr)E_{n,j,k} \\
    & & +\frac{E_{n,j,k}}{16\gamma^4}\Bigl(b_{d,\gamma} + (D^2_{n,j,k})^2 + e_{d,\gamma} D^2_{n,j,k}\Bigr)\Biggr].
\end{eqnarray*}
Also the class of tests based on ${\rm DEH}^*_{n,\gamma}$ is `closed at the boundaries $\gamma \to 0$ and $\gamma \to \infty$' since, elementwise on the underlying probability space, we have
\begin{equation}\label{limdehstar}
\lim_{\gamma \rightarrow 0} \left[\left(\frac{\gamma}{\pi}\right)^{d/2}\! {\rm DEH}^*_{n,\gamma} \! - \!  \frac{d(d\! +\! 2)}{4\gamma^2}\right] = b_{n,d}^{(2)}  \! - \! d^2, \ \
\lim_{\gamma \rightarrow\infty} \frac{2\gamma^{d/2+1}}{n\pi^{d/2}}{\rm DEH}^*_{n,\gamma} = \widetilde{b}_{n,d}^{(1)},
\end{equation}
where $b_{n,d}^{(2)}$ and $\widetilde{b}_{n,d}^{(1)}$ are given in \eqref{defskewnkurtma} and \eqref{defskmorosz}, respectively.
Under $H_0$, we have ${\rm DEH}^*_{n,\gamma} \vertk {\rm DEH}^*_{\infty,\gamma} := \int {\cal S}^2(t) \widetilde{w}_\gamma (t) \, {\rm d}t$, where ${\cal S}$ is some
centred Gaussian random element of ${\rm L}^2(\RR^d,{\cal B}^d,\widetilde{w}_\gamma(t){\rm d}t)$.
\cite{deh19b} also obtain a closed-form expression for $\BE [{\rm DEH}^*_{\infty,\gamma}]$.

If $X$ has a standardized alternative distribution satisfying $\BE \|X\|^4 < \infty$, we have
\[
\frac{{\rm DEH}^*_{n,\gamma}}{n} \fsk D^*_\gamma := \int |-\Delta\Psi^+(t) + (\|t\|^2 - d)\Psi^+(t)|^2 \widetilde{w}_\gamma(t) \, {\rm d}t,
\]
where $\Psi^+(t) = \mathbb{E}[\cos(t^\top X)] + \mathbb{E}[\sin(t^\top X)]$. Hence, $D^*_\gamma$ is the measure of distance from
$H_0$ associated with ${\rm DEH}^*_{n\gamma}$. Interestingly, under the stronger condition $\BE \|X\|^6 < \infty$, we
have
\[
\lim_{\gamma\rightarrow\infty}  \frac{2\gamma^{d/2+1}}{\pi^{d/2}}D^*_\gamma = \left\|\mathbb{E}\left(\|X\|^2X\right)\right\|^2.
\]
Since the right hand side is population skewness in the sense of \cite{morosz93} (see Section \ref{secskurt}), this result complements the
second limit in \eqref{limdehstar}. \cite{deh19b} also show that, under a fixed alternative distribution satisfying $\BE\|X\|^4 < \infty$,
$\sqrt{n}\big{(}{\rm DEH}^*_{n,\gamma}/n- D^*_\gamma\big{)})$ has a centred limit normal distribution with a variance that can be consistently estimated
from $X_1,\ldots,X_n$.

\section{The Henze--Zirkler test}\label{sechztest}
\cite{hezi90} observed that the BHEP-statistic defined in \eqref{defBHEP} may be written in the form
\[%\begin{equation}\label{reprBHEPden}
{\rm BHEP}_{n,\beta} = (2\pi)^{d/2} \beta^{-d} \int_{\R^d} \left(g_{n,\beta}(x) - \frac{1}{(2\pi \tau^2)^{d/2}} \exp\left(- \frac{\|x\|^2}{2 \tau^2} \right) \right)^2 {\rm d}x,
\]%\end{equation}
where $\tau^2 = (2\beta^2+1)/(2\beta^2)$, and
\[
g_{n,\beta}(x) = \frac{1}{nh^d} \sum_{j=1}^n \frac{1}{(2\pi)^{d/2}} \exp \left(- \frac{\|x- Y_{n,j}\|^2}{2h^2} \right),
\]
where $h^2= 1/(2\beta^2)$. The function $g_{n,\beta}$ is a nonparametric kernel density estimator with Gaussian kernel $w_1$  (recall $w_\beta$ from
\eqref{weightfunction}) and bandwidth $h$, applied to $Y_{n,1}, \ldots, Y_{n,n}$. A choice of the bandwidth $h$ in oder to minimize  the
mean integrated square error when estimating $w_1$ yields $h= h_n = (4/(2d+1)n)^{-1/(d+4)}$ and thus $\beta = \beta_n$, where
\begin{equation}\label{betahz}
\beta_n = 2^{-1/2} ((2d+1)n/4)^{1/(d+4)}.
\end{equation}
The Henze--Zirkler test statistic is given by ${\rm HZ}_n = {\rm BHEP}_{n,\beta_{n}}$. Apparently unaware of the work of \cite{hezi90}, \cite{bf93} proposed a test statistic ${\rm BF}_n$ that turned out to satisfy ${\rm BF}_n = \beta_n^d (2 \pi)^{d/2} {\rm BHEP}_{n,\beta_n}$ (see Section 7 of \cite{he02}.
Thus, ${\rm BF}_n$ is equivalent to a BHEP-statistic with a smoothing parameter that depends on $n$. \cite{gu00} proved that
\begin{equation}\label{bfh0asy}
\frac{nh^d2^d\pi^{d/2}{\rm BF}_n -1}{2^{1/2-d/4}h^{d/2}} \vertk {\rm N}(0,1)
\end{equation}
as $n \to \infty$ under $H_0$. Under a fixed standardized alternative distribution with density $f$, \cite{gu00} showed that
\begin{equation}\label{bfasymptotikalt}
\frac{\sqrt{n}}{2} \left({\rm BF}_n - \frac{1}{nh_n^d2^d\pi^{d/2}} - C(f,h_n) \right) \vertk {\rm N}(0,\sigma^2(f))
\end{equation}
for constants $\sigma^2(f)$ and $C(f,h_n)$, where $\lim_{n \to \infty} C(f,h_n) = \int (f(x)-w_1(x))^2 {\rm d}x$. In view of $nh_n^d \to \infty$,
\eqref{bfasymptotikalt} entails ${\rm BF}_n \stk \int (f(x)-w_1(x))^2 \, {\rm d}x$
under $f$. Hence, the test of $H_0$ based on ${\rm BF}_n$ (or ${\rm HZ}_n$) is consistent against general alternatives.
However,  since \eqref{bfh0asy} remains true under contiguous alternatives that approach $H_0$ at the rate $n^{-1/2}$,
the Henze--Zirkler (Bowman--Foster) test is not able to detect such alternatives, see also \cite{tenr07} for more general results
on Bickel--Rosenblatt-type statistics.

\section{The energy test}\label{secenergy}
For nearly 20 years now, the energy test has emerged as a strong genuine test for multivariate normality.
It is based on the notion of {\em energy distance} between multivariate distributions.
The naming {\em energy} stems from a close analogy with Newton's gravitational potential energy,  see, e.g., \cite{szri13}.
Besides goodness-of-fit testing, the concept of energy distance has found applications in many other fields, such as
testing for equality of distributions, nonparametric extensions of analysis of variance,
clustering, or testing for independence via distance covariance and distance correlation, see e.g., \cite{szri16}.

If $X$ and $Y$ are independent
random vectors with distributions $\PP^X$ and $\PP^Y$,and $X'$ and $Y'$ denote independent
copies of $X$ and $Y$, respectively, then the squared energy distance between $\PP^X$ and $\PP^Y$ is defined as
\[
D^2(\PP^X,\PP^Y) := 2 \BE\|X-Y\| - \BE \|X-X'\| - \BE\|Y-Y\|,
\]
provided these expectations exist (which is tacitly assumed). The energy distance $D(\PP^X,\PP^Y)$ satisfies
all axioms of a metric. A proof of the fundamental inequality $D(\PP^X,\PP^Y) \ge 0$, with  equality if
and only if $\PP^X=\PP^Y$, follows from \cite{zkk92} or \cite{ma97}, see also  \cite{szri05} for a different proof related to a result of  \cite{mo01}.

The energy test statistic for testing $H_0$ is
\[
{\cal E}_n := n \left(\frac{2}{n}\sum_{j=1}^n \BE\|\widetilde{Y}_{n,j}-N_1\|- \BE\|N_1-N_2\| - \frac{1}{n^2}\sum_{j,k=1}^n \|\widetilde{Y}_{n,j}-\widetilde{Y}_{n,k}\| \right).
\]
Here, $\widetilde{Y}_{n,j} = \sqrt{n/(n-1)}Y_{n,j}$ with $Y_{n,j}$ given in \eqref{scaledres},
and $N_1$ and $N_2$ are independent random vectors with the normal distribution N$_d(0,\text{I}_d)$, which are
independent of $X_1,\ldots,X_n$. The first expectation is with respect to $N_1$.
  Notice that $\BE\|N_1-N_2\| = 2\Gamma((d+1)/2)/\Gamma(d/2)$, where $\Gamma(\cdot)$ is the gamma function.
Since, for $a \in \RR^d$, the distribution of $\|a-N_1\|^2$ does only depend on $\|a\|^2$, the statistic ${\cal E}_n$ is seen to be invariant.
The energy test for multivariate normality rejects $H_0$ for large values of ${\cal E}_n$. It is consistent against each fixed non-normal alternative,
see \cite{szri05}, and it is fully implemented in the {\em energy} package for R, see \cite{risz14}. To the authors' knowledge,
there are hitherto no results on the behavior of ${\cal E}_n$ with respect to contiguous alternatives to $H_0$. Since the intrinsic (quadratic) measure
of distance between an alternative distribution $\PP^X$ (which, because of invariance, may be taken as having zero mean and unit covariance matrix)
and the standard $d$-variate normal distribution N$_d(0,\text{I}_d)$ is given by $\Delta_E(\PP^X) := D^2(\PP^X,{\text{N}}_d(0,{\text{I}}_d))$, say, it would be interesting to see
whether $\sqrt{n}({\cal E}_n - \Delta_E(\PP^X))$ has a non-degenerate normal limit as $n \to \infty$, with a variance that can consistently be estimated
from the data $X_1,\ldots,X_n$. Such a result would pave the way for an asymptotic confidence interval for  $\Delta_E(\PP^X)$.

\section{The test of Pudelko}\label{secpudelko}
For a fixed $r>0$, \cite{pu05} suggested to reject $H_0$ for large values of  the weighted supremum distance
\[
{\rm PU}_{n,r} = \sqrt{n} \sup_{0< \|t\| \le r} \frac{|\Psi_n(t)- \Psi_0(t)|}{\|t\|},
\]
where $\Psi_n(t)$ is given in \eqref{defpsint}, and $\Psi_0(t) = \exp(-\|t\|/2)$. The test statistic ${\rm PU}_{n,r}$ is invariant,
since it is a function of the scaled residuals $Y_{n,1},\ldots,Y_{n,n}$ and rotation invariant. This statistic is similar in spirit
as the statistic studied by \cite{cs86}, which is $\sup_{\|t\|\le r}\big{|}|\Psi_n(t)|^2 - \Psi_0^2(t)\big{|}$. Under $H_0$,
${\rm PU}_{n,r}$ converges in distribution to $\sup_{0<\|t\|\le r} |{\cal P}(t)|/\|t\|$, where
${\cal P}(\cdot)$ is a centred Gaussian random element of the Banach space $C(B_r)$ of complex-valued continuous functions,
defined on $B_r:= \{x \in \RR^d: \|x \| \le r\}$, equipped with the
supremum norm $\|f\|_{C(B_r)} := \sup_{x \in B_r}|f(x)|$. \cite{pu05} also showed that the test is able to detect contiguous alternatives
that approach $H_0$ at the rate $n^{-1/2}$. The consistency of the test based on ${\rm PU}_{n,r}$ follows easily from
\cite{cs89}. A drawback of this test is its lack of feasibility, since one has to calculate the supremum of a function inside a
$d$-dimensional sphere.

\section{The test of Cox and Small}\label{seccoxsmall}
According to \cite{cs78}, a main objective of tests of $H_0$ is 'to see whether an estimated covariance matrix
 provides an adequate summary of the interrelationships among a set of variables', and that departure from
 multivariate normality 'is often the occurrence of appreciable nonlinearity of dependence'.
 To obtain an affine invariant test that assesses the degree of nonlinearity,
    they propose to find that pair of linear combinations of the original variables,
   such that one has maximum curvature in its regression on the other.
  The population functional which underlies the test of Cox and Small is $T_{CS}(\PP^X)=\max_{b\in\mathcal{S}^{d-1}}\eta^2(b)$, where
\begin{equation*}
\eta^2(b)=\frac{\left\|\BE\left(X(b^\top X)^2\right)\right\|^2-\left(\BE\left(b^\top X\right)^3\right)^2}{\BE\left(b^\top X\right)^4-1-\left(\BE\left(b^\top X\right)^3\right)^2},
\end{equation*}
see \cite{cs78}, p. 268. The test statistic is
$
T_{n,CS}=\max_{b\in\mathcal{S}^{d-1}}\eta_n^2(b)$,
where
\begin{equation*}
\eta_n^2(b)=\frac{\left\|n^{-1}\sum_{j=1}^nY_{n,j}(b^\top Y_{n,j})^2\right\|^2-\left(n^{-1}\sum_{j=1}^n(b^\top Y_{n,j})^3\right)^2}{n^{-1}\sum_{j=1}^n(b^\top Y_{n,j})^4-1-\left(n^{-1}\sum_{j=1}^n(b^\top Y_{n,j})^3\right)^2}
\end{equation*}
is the empirical counterpart of $\eta^2(b)$.
Rejection of $H_0$ will be for large values of $T_{n,CS}$.
The statistic $T_{n,CS}$ is affine invariant, since it is both a function of $Y_{n,1},\ldots,Y_{n,n}$ and rotation invariant. Notice that
the functional  $T_{CS}$ vanishes on the set ${\cal N}_d$, but $T_{CS}(\PP^X) =0$ does not necessarily imply that $\PP^X \in {\cal N}_d$.
Some missing distributional properties of the statistic $T_{n,CS}$ were provided by \cite{eb12}. If $\PP^X$ is elliptically symmetric and satisfies $\BE \|X\|^6 < \infty$,
then
\begin{equation*}
nT_{n,CS}\vertk \frac{d(d+2)}{3m_4-d(d+2)}\max_{b\in\mathcal{S}^{d-1}}W(b)^\top B W(b),
\end{equation*}
where $m_4=E\|X\|^4$, $B$ is the $(d+1)\times(d+1)$-matrix $\mbox{diag}(1,\ldots,1,-1)$,
 and $W(\cdot)$ is a centred $(d+1)$-variate Gaussian process in $C(\mathcal{S}^{d-1},\R^{d+1})$, the space of continuous functions
 from $\mathcal{S}^{d-1}$ to $\R^{d+1}$ (see Theorem 2.4 of \cite{eb12}, where the covariance matrix kernel of $W$ is given explicitly).
  As a consequence, the test of Cox and Small is not able to detect such elliptical alternatives to normality.
  Next, writing $\mu(b)=\BE((b^\top  X)^2(X,(b^\top  X))^\top)$, we have
\[%\begin{equation}\label{eq:CSconvP}
T_{n,CS}\stk \max_{b\in\mathcal{S}^{d-1}}\frac{\mu(b)^\top B\mu(b)}{\BE(b^\top  X)^4-1-(\BE(b^\top  X)^3)^2}
\]%\end{equation}
if $\BE \|X\|^6 < \infty$. Thus, the test based on
   $T_{n,CS}$ is consistent against each alternative distribution for which the above stochastic limit $\delta(\PP^X)$ (say) is positive.
 \cite{eb12} also provides the limit distribution of $T_{n,CS}$ under contiguous alternatives to $H_0$, but it
 is still  an open problem whether $\sqrt{n}(T_{n,CS}- \delta(\PP^X))$ has
 a nondegenerate limit distribution as $n \to \infty$. From a practical point of view, the test of Cox and Small
 has the drawback that finding the maximum of $\eta_n^2(b)$ over $b \in \mathcal{S}^{d-1}$ is a computationally extensive task.

\section{The test of Manzotti and Quiroz}\label{secquiroz}
 \cite{maqu01} propose to test $H_0$ by means of averages over the standardized sample of multivariate spherical harmonics,
 radial functions and their products. For $k\in\mathbb{N}$ let $f_1,\ldots,f_k:\R^d \to \R$, such that $\BE f_j^2(X) <\infty$ if
 $X\edist  \mbox{N}_d(0,{\rm I}_d)$, $j=1,\ldots,k$. Let $V=(v_{ij})$ be the ($k\times k$)-matrix with entries
\begin{equation*}
v_{ij}=\BE[f_i(X)f_j(X)]-\BE f_i(X) \,  \BE f_j(X),\quad X\edist \mbox{N}_d(0,I_d),\end{equation*}
where $V$ is assumed to be invertible. For $\mathbf{f}=(f_1,\ldots,f_k)^\top$, let
\begin{equation*}
\nu_n(f_j)=\frac1{\sqrt{n}}\sum_{\ell=1}^n \big{\{} f_j(Y_{n,\ell})-\BE f_j(X) \big{\}}  \quad\mbox{and}\quad \nu_n(\mathbf{f})=(\nu_n(f_1),\ldots,\nu_n(f_k))^\top .
\end{equation*}

The general type of test statistic of \cite{maqu01} is  the quadratic form
\begin{equation*}
T_{n,MQ}(\mathbf{f})=\nu_n(\mathbf{f})^\top V^{-1}\nu_n(\mathbf{f}).
\end{equation*}
To be more specific, let  $\mathcal{H}_j$, $j \ge 0$,  be the set of spherical harmonics of degree $j$
 in the orthonormal basis of spherical harmonics in $d$ dimensions with respect to the uniform measure on $\mathcal{S}^{d-1}$, and put $\mathcal{G}_j=\bigcup_{i=0}^j\mathcal{H}_i$.
 The number of linear independent spherical harmonics of degree $j$ in dimension $d$ is ${d+j-1 \choose j} - {d+j-3 \choose j-2}$.
  A suitable orthonormal basis can be found using Theorem 5.25 in \cite{ABW:2001} or \cite{maqu01}, see also \cite{G:1996} or \cite{M:1998} for  details on spherical harmonics.
   \cite{maqu01} suggest two different choices for $\mathbf{f}$. Putting $r_j(x)=\|x\|^j$, $x \in \R^d$, and $u(x)=x/\|x\|$, $x \neq 0$,
 the first statistic $T_{n,MQ}(\mathbf{f}_1)$ uses $f_j$ of the form $g\circ u$ for $g\in\mathcal{G}_4\setminus \mathcal{H}_0$, giving a total of $k={d+3 \choose 4} - {d+2 \choose 3}-1$ functions.
 Due to orthonormality we have $V={\rm I}_k$, and since no radial functions are considered, $T_{n,MQ}(\mathbf{f}_1)$ only tests for aspects of spherical symmetry.
The second statistic $T_{n,MQ}(\mathbf{f}_2)$ uses the functions $r_1$ and $r_3 (g \circ u)$, where $g \in \mathcal{G}_2$, which comprise a totality of
$k={d+1 \choose 2} +d + 1$ functions.

Both statistics are affine invariant, and \cite{maqu01} derive their limit null distributions, which are sums of weighted  independent $\chi^2_1$ random variables.
Although the authors do not deal with the question of consistency of their tests, it is easily seen that, under an alternative distribution $\PP^X$ (which, in view of invariance, is assumed
to be standardized), and suitable conditions on $f_1,\ldots,f_k$,
we have
\[
\frac{1}{n} T_{n,MQ} \stk \delta(f)^\top V^{-1}   \delta(f)
\]
as $n \to \infty$, where $\delta(f) = (\BE f_1(X) - \BE_0 f_1, \ldots,\BE f_k(X) - \BE_0 f_k)^\top$,  and
$\BE_0 f_j$ is the expectation $\BE f_j(N)$, where $N \edist {\rm N}_d(0,{\rm I}_d)$.
Since there are non-normal distributions for which the  above (non-negative) stochastic limit vanishes,
the tests of \cite{maqu01} are not consistent against general alternatives.
To the best of our knowledge, there are no further asymptotic properties of $T_{n,MQ}$ under alternatives to $H_0$.

\section{Tests based on skewness and kurtosis}\label{secskurt}
A still very popular group of tests for $H_0$ employ  measures of multivariate skewness and kurtosis. The popularity of
these tests stems from the widespread belief that, in case of rejection of $H_0$, there is some evidence regarding the kind of departure
from normality of the underlying distribution. The then state of the art regarding this group of  tests has been reviewed in \cite{he02}, but for the sake of
completeness, we revisit the most important facts. The classical invariant  measures of multivariate sample skewness and
kurtosis  due to  \cite{ma70} are defined by
\begin{equation}\label{defskewnkurtma}
b_{n,d}^{(1)} = \frac{1}{n^2} \sum_{j,k=1}^n \left(Y_{n,j}^\top Y_{n,k}\right)^3, \qquad b_{n,d}^{(2)} = \frac{1}{n} \sum_{j=1}^n \|Y_{n,j}\|^4,
\end{equation}
respectively. The functional (population counterpart) corresponding to $b_{n,d}^{(1)}$ is $\beta_d^{(1)} = \beta_d^{(1)}(\PP^X) = \BE (X_1^\top X_2)^3$,
where $X$ is standardized, $X_1,X_2$ are i.i.d. copies of $X$, and $\BE \|X\|^6 < \infty$. The functional accompanying kurtosis
is $\beta_d^{(2)} = \beta_d^{(2)}(\PP^X) = \BE \|X\|^4$, where, like above, $\BE(X) =0$ and $\BE(XX^\top) = {\rm I}_d.$ When used as statistics
to test $H_0$,  $b_{n,d}^{(1)}$ has an upper rejection region, whereas the test based on $b_{n,d}^{(2)}$ is two-sided.
If the distribution of $X$ is elliptically symmetric, we have
\begin{equation}\label{limitmardiaskew}
n b_{n,d}^{(1)} \vertk \alpha_1 \chi_d^2 + \alpha_2\chi_{d(d-1)(d+4)}^2,
\end{equation}
where
\[
\alpha_1 = \frac{3}{d}\bigg{[} \frac{\BE\|X\|^6}{d+2} - 2 \BE\|X\|^4 + d(d+2)\bigg{]}, \qquad \alpha_2 = \frac{6\BE\|X\|^6}{d(d+2)(d+4)},
\]
where $\chi_d^2$, $\chi_{d(d-1)(d+4)}^2$ are independent $\chi^2$-variables with $d$ and $d(d-1)(d+4)$ degrees of freedom, respectively,
see \cite{bahe92}, and \cite{kl02}. Notice that  $\alpha_1=\alpha_2 =6$ under $H_0$, whence $n b_{n,d}^{(1)} \vertk 6 \chi^2_{d(d+1)(d+2)/6}$
under normality, see \cite{ma70}.  From \eqref{limitmardiaskew}, it follows that the test of $H_0$ based on $b_{n,d}^{(1)}$ is not consistent against
spherically symmetric alternatives satisfying $\BE\|X\|^6 < \infty$. If $\beta_d^{(1)} > 0$, then $\sqrt{n}(b_{n,d}^{(1)} - \beta_d^{(1)})$ has a centred non-degenerate
limit normal distribution as $n \to \infty$, see Theorem 3.2 of \cite{bahe92}.
The skewness functional $\beta_d^{(1)}(\cdot)$ does not characterize the class ${\cal N}_d$ of normal distributions since,
although $\beta_d^{(1)}(\cdot)$ vanishes on ${\cal N}_d$, there are (notably elliptically {\em symmetric}) non-normal
distributions that share this property. Since the critical value of $b_{n,d}^{(1)}$ as a test statistic for assessing
multivariate normality is computed under the very assumption of {\em normality},
the inclination to impute supposedly diagnostic properties to $b_{n,d}^{(1)}$ in case of rejection of $H_0$ in the sense
that 'there is evidence that the underlying distribution is skewed' is not justified, at least not in terms of statistical significance.
In fact, the limit
distribution of $n b_{n,d}^{(1)}$ under certain classes of elliptically symmetric distributions is
stochastically much larger than the limit null distribution of $n b_{n,d} ^{(1)}$ (see \cite{bahe92}), and so rejection of $H_0$ based on $b_{n,d}^{(1)}$ may be due to an underlying long-tailed
elliptically {\em symmetric} distribution.

Regarding kurtosis, we have
$\sqrt{n}(b_{n,d}^{(2)} - \beta_d^{(2)}) \vertk {\rm N}(0,\sigma^2)$ as $n \to \infty$, where $\sigma^2$ depends on
mixed moments of $X$ up to order 8, see \cite{he94}. Under $H_0$, we have $\beta_d^{(2)} = d(d+2)$ and
$\sigma^2 = 8d(d+2)$, and the limit distribution was already obtained by \cite{ma70}, see also \cite{kl02} for the case
that $\PP^X$ is elliptically symmetric. It follows that, under the
condition $\BE \|X\|^8 < \infty$,  Mardia's kurtosis test for normality is consistent if and only if $\beta_d^{(2)} \neq d(d+2)$.
The critical remarks made above on alleged diagnostic capabilies of tests for $H_0$ based on measures of skewness apply
mutatis mutandis to a test for normality  based on $b_{n,d}^{(2)}$ or any other measure of multivariate kurtosis.

Among the many measures of multivariate skewness, we highlight skewness in the sense of \cite{morosz93}, because
it emerges in connection with several weighted $L^2$-statistics for testing $H_0$. This measure is defined by
\begin{equation}\label{defskmorosz}
\widetilde{b}_{n,d}^{(1)} := \frac{1}{n^2} \sum_{j,k=1}^n \|Y_{n,j}\|^2 \|Y_{n,k}\|^2 Y_{n,j}^\top Y_{n,k}.
\end{equation}
The corresponding functional (population counterpart) is $\widetilde{\beta}_d^{(1)}= \big{\|} \BE (\|X\|^2 X)\big{\|}^2$,
where $X$ is assumed to be standardized and $\BE \|X\|^6 < \infty$. Limit distributions for $\widetilde{b}_{n,d}^{(1)}$
have been obtained by \cite{he97} both for the case that $\PP^X$ is elliptically symmetric (which implies  $\widetilde{\beta}_d^{(1)}=0$)
and the case that $\widetilde{\beta}_d^{(1)}>0$, see also \cite{kl02}. A further measure of multivariate skewness that has  been reviewed in \cite{he02} is skewness in the sense of \cite{maaf73}, which is defined
as
\[
b_{n,d,M}^{(1)} = \max_{u \in \mathcal{S}^{d-1}} \frac{\big{\{}n^{-1}\sum_{j=1}^n (u^\top X_j- u^\top \overline{X}_n)^3 \big{\}}^2}{(u^\top S_n u)^3}.
\]
General limit distribution theory for $b_{n,d,M}^{(1)}$ is given in \cite{bahe91}. As for further measures of multivariate kurtosis, we mention the measure
\[
\widetilde{b}_{n,d}^{(2)} = \frac{1}{n^2} \sum_{j,k=1}^n \left(Y_{n,j}^\top Y_{n,k}\right)^4,
\]
introduced by \cite{ko89}. The corresponding functional is $\widetilde{\beta}_d^{(2)} = \BE(X_1^\top X_2)^4$, where
$X_1,X_2$ are i.i.d. copies of the standardized vector $X$, and $\BE\|X\|^8 < \infty$. General asymptotic distribution theory for
$\widetilde{b}_{n,d}^{(2)}$ is provided by \cite{he94b} and \cite{kl02}. \cite{he02} also reviewed kurtosis in the sense of
\cite{maaf73}, which is defined as
\[
b_{n,d,M}^{(2)} = \max_{u \in \mathcal{S}^{d-1}} \frac{n^{-1}\sum_{j=1}^n (u^\top X_j- u^\top \overline{X}_n)^4}{(u^\top S_n u)^2}.
\]
Limit distribution theory for $b_{n,d,M}^{(2)}$ has been obtained by \cite{bahe91} and \cite{naito98}.

Since the review \cite{he02}, there have been the following suggestions to test $H_0$ by means of measures of
multivariate skewness and kurtosis (which, however, do not lead to consistent tests and share the drawback stated at the beginning of this section):
 \cite{kto07} consider invariant tests of multivariate normality that are based on the Mahalanobis distance between two multivariate location vector estimates (as a measure of skewness) and on the (matrix) distance between two scatter matrix estimates (as a measure of kurtosis). Special choices of these estimates yield generalizations of Mardia's skewness an kurtosis. The authors obtain asymptotic distribution theory of their test statistics both under normality and
certain contiguous alternatives to $H_0$, and they compare the limiting Pitman efficiencies to those of Mardia's tests based on
 $b_{n,d}^{(1)}$ and $b_{n,d}^{(2)}$. \cite{doha08} propose a non-invariant test based on skewness and kurtosis. \cite{ehhs} consider a transformation of Mardia's kurtosis statistic, with the aim of improving the finite-sample approximation with respect to a normal limit distribution.

\section{Miscellaneous results}\label{secmisc}
\cite{arco07} proposed two invariant test statistics that are based on the following characterizations, see, e.g., \cite{cr36}.
Let $m \ge 2$ be a fixed integer, and let $X_1,\ldots,X_m$ be i.i.d. $d$-dimensional vectors satisfying $\BE(X_1) =0$ and $\BE (X_1X_1^\top) = {\rm I}_d$.
Then $m^{-1/2}\sum_{j=1}^m X_j \edist {\rm N}_d(0,{\rm I}_d)$ if and only if $X_1 \edist {\rm N}_d(0,{\rm I}_d)$. Furthermore,
$m^{-1/2}\sum_{j=1}^m X_j \edist X_1$ if and only if $X_1 \edist {\rm N}_d(0,{\rm I}_d)$. A statistic that corresponds to the first characterization is
$
\widehat{D}_{n,m} = \textstyle{\int} \big{|} \widehat{\Psi}_{n,m}(t) - \Psi_0(t)\big{|}^2 w_\beta(t) \, {\rm d}t
$,
where
$
\widehat{\Psi}_{n,m}(t) = $
$n!^{-1}(n-m)! \textstyle{\sum_{\neq}} \exp\big{(}{\rm i} t^\top m^{-1/2} \sum_{p=1}^m Y_{n,j_p} \big{)}$,
and $\Sigma_{\neq}$ means summation over all $j_1,\ldots,j_m \in \{1,\ldots,n\}$ such that $j_p \neq j_q$ if $p \neq q$.
Notice that this approach is a generalization of the BHEP-statistic given in \eqref{defBHEP}.
The statistic which is tailored to   the second characterization is
$
\widehat{E}_{n,m} = \textstyle{\int} \big{|} \Psi_{n,m}(1) - \Psi_{n,1}(t)\big{|}^2 w_\beta(t) \, {\rm d}t$.
Both statistics have representations in form of multiple sums.
By using the theory of $U$-statistics with estimated parameters, \cite{arco07} derives almost sure limits
of $\widehat{D}_{n,m}$ and $\widehat{E}_{n,m}$ as well as the limit distributions
of $n\widehat{D}_{n,m}$ and $n\widehat{E}_{n,m}$ under $H_0$. Some very limited simulations, performed for $n \le 15$ and $d=2$, indicate that
the power of these tests is comparable to that of the BHEP-test. However, the computational burden involved increases rapidly with $m$. 

Without providing any distribution theory, \cite{hhr02} suggest an invariant two-stage test procedure for testing $H_0$.
This procedure combines a modified correlation coefficient related to a Q-Q-plot of the ordered values of
$\|Y_{n,j}\|^2$, $j=1,\ldots,n$, against ordered quantiles of the $\chi^2_d$-distribution, and a test based on Mardia's
non-negative invariant measure of skewness $b_{n,d}^{(1)}$ given in \eqref{defskewnkurtma}. \cite{lpy04} deal with Q-Q-plots based on functions of $(j(j+1))^{-1/2}(X_1+ \ldots + X_j - jX_{j+1})$, $j=1,\ldots,n-1$, and hence recommend procedures that are not even invariant with respect to permutations of $X_1,\ldots,X_n$. The latter objection also holds for the procedure suggested by \cite{libe99}. \cite{tftw05} extend the projection procedure of \cite{lili00} to test for multivariate normality with incomplete longitudinal data with small sample size, including cases when the sample size $n$ is smaller than $d$. \cite{hata08} correct an inaccuracy of the (non-invariant) test of \cite{srhu87}, and \cite{maru07} derives approximations of expectations
and variances related to that test under alternative distributions. Without providing any theoretical results, \cite{hata08} aim at transforming two graphical methods for assessing $H_0$ into formal statistical tests. A variant of this approach was considered by \cite{maok18}. \cite{cdf10} suggest to perform a chi-quare test based on $\|Y_{n,1}\|^2, \ldots, \|Y_{n,n}\|^2$ (see also \cite{most81}),  and \cite{bmpz13} extend this approach to include more general power divergence type of test statistics. \cite{maok19} consider $\ell_1$- and $\ell_2$-type measures of deviation between $\|Y_{n,j}\|^2$ and corresponding approximate expected order
statistics of a $\chi^2_d$-distribution (for tests based on $\|Y_{n,1}\|^2, \ldots, \|Y_{n,n}\|^2$, see also Section 5.2 of \cite{he02}).
\cite{vpmv16} compare several test statistics that, for fixed $r \ge 2$,  are quadratic forms in the vector $(V_{n,1}, \ldots, V_{n,r})^\top$.
Here, $V_{n,j} = (N_{n,j}-n/r)/(\sqrt{n/r})$, $N_{n,j} = \sum_{k=1}^n {\bf 1}\{c_{j-1} < \|Y_{n,k} \|^2 \le c_j\}$, and $0< c_1 < \ldots < c_{r-1} < c_r = \infty$,
where $c_j$ is the ($j/r$)-quantile of the $\chi^2_d$-distribution, $j=1,\ldots,r-1$. \cite{joen11} investigates the finite-sample performance of of the Jarque--Bera test for $H_0$ in order to improve the size of the test. \cite{khp14} improve upon multivariate Jarque--Bera type tests by means of transformations.
Simulations show that such transformatinos essentially improves test accuracy when $d$ is close to $n$.
\cite{kim16} generalizes the univariate Jarque--Bera test and its modifications
to the multivariate versions using an orthogonalization of data and compares
it with competitors in a simulation study. \cite{kipa18} propose a non-invariant test based on univariate Anderson--Darling type statistics that are averaged out over the $d$ coordinates. \cite{vage09} suggest a non-invariant test that is based on the average of Shapiro--Wilk statistics, applied to each of the components of
$Y_{n,1},\ldots,Y_{n,n}$. By using an idea of \cite{frla06}, \cite{tenr11} proposes an invariant consistent multiple test procedure that
combines Mardia's measures of skewness and kurtosis and two members of the family of BHEP tests. The combined procedure
rejects $H_0$ if one of the statistics is larger than its $(1-u_{n,\alpha})$-quantile under $H_0$, where
$u_{n,\alpha}$ is calibrated so that the combined test has a desired level of significance  $\alpha$.
In the same spirit, \cite{tenr17} combines two BHEP-tests and the 'extreme' BHEP-tests, the statistics of which are
given by the right hand sides of \eqref{bhepbeta0} and  \eqref{bhepbetainf}. \cite{masz10} consider the problem of testing $H_0$ against some alternatives that are invariant with respect to a subgroup of the full group of affine transformations and obtain
approximations to the most powerful invariant tests. Special emphasis is given to
exponential and uniform alternatives in the case $d=2$, whereas the case $d \ge 3$ is
only sketched. In the spirit of projection pursuit tests (see Section 8.1 of \cite{he02}),
which are  based on Roy's union-intersection principle (\cite{ro53}), \cite{zhsh14}
propose a non-invariant test that combines the Shapiro--Wilk test and Mardia's kurtotis test.
In the same spirit, \cite{wahw11} suggest a statistic that considers solely the Shapiro--Wilk statistic.

\cite{wang14} provides a MATLAB package for testing $H_0$, which is  implemented as an interactive and graphical tool.
The package comprises 12 different tests, among which are the energy test, the Henze--Zirkler test, and the tests
based on Mardia's skewness and kurtosis. \cite{thu14} proposes six invariant tests for $H_0$, the common basis of which are characterizations of independence
of sample moments of the multivariate normal distribution.

\section{Comparative simulation studies}\label{secsimstud}

\subsection{Available simulation studies}
\cite{memu05} perform an extensive simulation study with 13 tests for multivariate normality.
From this study, they conclude that 'if one is going to rely on one and only one procedure, the Henze--Zirkler test is
recommended. This recommendation is based on the relative ease of use (the test statistic has
an approximately lognormal asymptotic distribution), good Monte Carlo simulation results,
and mathematically proven consistency against all alternatives'. \cite{fsbn07} compare four tests of multivariate normality and conclude: 'The results of our simulation suggest that, relative to the other two tests considered, the Henze and Zirkler test generally possesses good power across the alternative distributions
investigated, in particular for $n \ge 75$'. \cite{hesk18} compare four test of $H_0$ that are based on a combination of measures of multivariate skewness and kurtosis, and the Henze--Zirkler test. They concluded that 'the Henze--Zirkler test best preserves the nominal significance level', and that
'for the number of traits and sample sizes considered, it is not possible to indicate the most
powerful test for all kinds of alternative distributions considered in the paper'. \cite{jovo14} investigate 15 tests of $H_0$, all of which freely available as R-functions. They find that some tests are unreliable and should either be corrected or removed, or their deficits should be commented
upon in the documentation by the package maintainer. Moreover, they summarize:
'On the question of whether or not multivariate tests offer an advantage over simply testing each
marginal distribution with a univariate test, the answer is a resounding yes. Not only are some
multivariate tests able to detect deviations from normality that are not reflected in the marginals
of the distribution, but these tests are also, in part, more powerful for distributions that do display
the deviations in the marginals'.

\subsection{New simulation study}
 This subsection compares the finite-sample power performance of the tests presented in this survey by means of a Monte Carlo simulation study. All simulations are performed using the statistical computing environment \textsf{R}, see \cite{r20}. The tests were implemented in the accompanying {\tt R} package {\tt mnt}, see \cite{be20}.

 We consider the sample sizes $n=20$, $n=50$ and $n=100$, the dimensions $d=2$, $d=3$ and $d=5$, and the nominal level of significance is set to $0.05$.
  Throughout, critical values for the tests have been simulated with $100~000$ replications under $H_0$, see Table \ref{tab:cv}.
  Note that, in order to ease the comparison with  the original articles, we state the empirical quantiles of
  $\left(16\gamma^{2+d/2}/\pi^{d/2}\right){\rm HV}_{n,\gamma}$, $\pi^{-d/2}{\rm HJ}_{n,\gamma}$,
  $(\gamma/\pi)^{d/2}{\rm HJM}_{n,\gamma}$, $(\gamma/\pi)^{d/2}d^{-2}{\rm DEH}_{n,\gamma}$, and
  $(\gamma/\pi)^{d/2}d^{-2}{\rm DEH}^*_{n,\gamma}$ and chose whenever available the tuning parameter $\gamma$ according to the suggestions of the authors, respectively.
  For the sake of readability,  we subduct the index $n$ for all tests in the tables.
  The values of Table \ref{tab:cv} are also reported in package {\tt mnt} in the data frame {\tt Quantile095} for easy access.
  Each entry in a table that refers to empirical rejection rates as estimates of the power of the test is based on $10~000$ replications,
  with the exception of the HJM test, where $1~000$ replications have been considered, due to the heavy computation time of the procedure.
\begin{table}[t]
\small
\setlength{\tabcolsep}{1mm}
\centering
\begin{tabular}{cc|ccccccccc}
  $d$ & $n$ & $b^{(1)}$ & $b^{(2)}$ & $b_{M}^{(1)}$ & $\widetilde{b}^{(2)}$ & $\widetilde{b}^{(1)}$ & $b_{M}^{(2)}$ & BHEP$_{1}$ & HZ & HV$_{5}$\\[1mm]
  \cline{1-11}
 \multirow{3}{*}{$2$} & 20 & 2.38 & 9.44 & 1.82 & 40.90 & 1.77 & 5.47 & 0.54 & 0.73 & 250 \\
 & 50 & 1.09 & 9.44 & 0.84 & 37.28 & 0.87 & 4.94 & 0.55 & 0.88 & 358 \\
 & 100 & 0.56 & 9.17 & 0.43 & 33.40 & 0.46 & 4.42 & 0.56 & 0.97 & 397\\[1mm]
 \multirow{3}{*}{$3$} & 20 & 4.63 & 16.37 & 2.81 & 75.36 & 2.68 & 6.68 & 0.67 & 0.82 & 545\\
  & 50 & 2.11 & 16.73 & 1.25 & 67.38 & 1.39 & 5.81 & 0.68 & 0.92 & 823\\
  & 100 & 1.09 & 16.49 & 0.63 & 60.05 & 0.74 & 5.02 & 0.68 & 0.98 & 936\\[1mm]
 \multirow{3}{*}{$5$} & 20 & 12.57 & 35.35 & 4.38 & 191 & 4.55 & 8.37 & 0.84 & 0.91 & 1750\\
  & 50 & 5.77 & 37.01 & 1.96 & 163 & 2.61 & 7.16 & 0.85 & 0.96 & 2993\\
  & 100 & 2.96 & 36.94 & 0.94 & 140 & 1.44 & 5.92 & 0.85 & 0.99 & 3530\\[2mm]

  &  & HJ$_{1.5}$ & HJM$_{1.5}$ & DEH$_{0.25}$ & DEH$^*_{0.5}$ & ${\cal E}$ & $T_{MQ}(\mathbf{f}_1)$ & $T_{MQ}(\mathbf{f}_2)$ & $T_{CS}$ & PU$_{2}$ \\[1mm]
  \cline{3-11}
 \multirow{3}{*}{$2$} & 20 & 12.31 & 2.89 & 1.92 & 3.42 & 0.93 & 11.26 & 3.71 & 0.38 & 1.02 \\
 & 50  & 42.83 & 3.40 & 1.98 & 3.50 & 0.96 & 11.17 & 4.39 & 0.16 & 1.02 \\
 & 100 & 80.57 & 3.62 & 1.96 & 3.53 & 0.97 & 11.27 & 4.65 & 0.08 & 1.03 \\[1mm]
 \multirow{3}{*}{$3$} & 20 & 32.74 & 6.67 & 1.60 & 2.44 & 1.04 & 29.13 & 3.91 & 0.59 & 1.18 \\
  & 50 &  148 & 9.32 & 1.66 & 2.52 & 1.07 & 29.09 & 4.67 & 0.26 & 1.19 \\
  & 100  & 335 & 9.71 & 1.65 & 2.53 & 1.07 & 29.05 & 4.99 & 0.13 & 1.20 \\[1mm]
 \multirow{3}{*}{$5$} & 20  & 127 & 25.64 & 1.36 & 1.79 & 1.23 & 115 & 4.52 & 0.82 & 1.33 \\
  & 50 & 1049 & 55.62 & 1.42 & 1.85 & 1.26 & 113 & 5.23 & 0.42 & 1.35 \\
  & 100 & 3117 & 72.65 & 1.42 & 1.86 & 1.28 & 113 & 5.61 & 0.22 & 1.36
\end{tabular}
\caption{Empirical $95\%$ quantiles of the test statistics under $H_0$ ($100~000$ replications)}\label{tab:cv}
\end{table}

We consider a total of 29 alternatives as well as a representative of the multivariate normal distribution. By NMix$(p,\mu,\Sigma)$
we denote the normal mixture distribution generated by
\begin{equation*}
	(1 - p) \, {\rm N}_d(0, {\rm I}_d) + p \, {\rm N}_d(\mu, \Sigma), \quad p \in (0, 1), \, \mu \in \R^d, \, \Sigma > 0,
\end{equation*}
where $\Sigma > 0$ stands for a positive definite matrix. In the notation of above, $\mu=3$ stands for a $d$-variate vector of 3's
and $\Sigma={\rm B}_d$ for a $(d \times d)$-matrix containing 1's on the main diagonal and 0.9's for each off-diagonal entry.
  We write  $t_\nu(0,{{\rm I}}_d)$ for the multivariate $t$-distribution with $\nu$ degrees of freedom, see \cite{GB09}.
 By DIST$^d(\vartheta)$ we denote the $d$-variate random vector generated by independently simulated components of the distribution DIST
 with parameter vector $\vartheta$, where DIST is taken to be the uniform distribution U, the lognormal distribution LN, the beta distribution B,
 as well as the Pearson Type II P$_{II}$ and Pearson Type VII distribution P$_{VII}$. For the latter distribution, we used the {\tt R} package {\tt PearsonDS}, see \cite{bk17}.
 The spherical symmetric distributions were simulated using the {\tt R} package {\tt distrEllipse}, see \cite{RKSC06}, and they
 are denoted by $\mathcal{S}^d(\mbox{DIST})$, where DIST stands for the distribution of the radii, which was chosen to be the exponential,
  the beta, the $\chi^2$-distribution and the lognormal distribution. With MAR$_d$(DIST) we denote ${\rm N}_d(0, {\rm I}_d)$-distributed random vectors,
 where the $d$th component is independently replaced by a random variable following the distribution DIST.
 Here, we chose the exponential, the $\chi^2$, student's $t$ and the gamma distribution. With NM$_d(\vartheta)$ we denote the normal mixture distributions generated by
\begin{equation*}
	0.5 \, {\rm N}_d(0, \Sigma_\vartheta) + 0.5 \, {\rm N}_d(0, \Sigma_{-\vartheta}),
\end{equation*}
where $\Sigma_\vartheta$ is a positive definite ($d\times d$)-matrix with 1's on the diagonal and the constant $\vartheta$ for each off diagonal entry.
 In this family of non-normal distributions each component follows a normal law. The symbol S$|\mbox{N}_d|$ stands for the distribution of $\pm|X|$,
 where $X\edist {\rm N}_d(0,I_d)$, the absolute value $|\cdot|$ is applied componentwise,  and $\pm$ assigns,  independently of each other and
 with equal probability 0.5,  a random sign to each component of $|X|$. Finally, we consider  the distribution ${\rm N}_d(\mu_d,\Sigma_{0.5})$,
 with $\mu_d=(1,2,\ldots,d)^\top$ and the same covariance structure as reported for the NM-alternatives, in order to show that all tests under consideration
 are invariant and indeed have a type I error equal to the significance level of $5\%$.

The results of the weighted $L^2$-type tests in Tables \ref{tab:d2} - \ref{tab:d5} are presented for the same tuning parameters as in Table \ref{tab:cv},
 and in order to keep the tables concise the values are omitted.

 First, we evaluate the results for $d=2$. A close look at Table \ref{tab:d2} reveals that, for the family of normal mixture distributions,
  the HZ-test and the PU-test perform best when the shifted standard normal distributions are mixed, whereas for different covariance matrices,
 the strongest procedure is HJM. The HJM-test performs also best throughout the multivariate $t$-distributions.
 For the independently simulated components, $T_{MQ}(\mathbf{f}_2)$ is strong, especially for marginal distributions with bounded support.
 Interestingly, each of the tests that are based on measures of skewness and kurtosis, as well as the HV- and the HJ-test,
 completely fail to detect these alternatives. For the Pearson-Type VII alternatives, HJM again has the strongest power,
while BHEP shows the strongest performance for LN$^2(0,0.5)$ and B$^2(1,2)$. The spherically symmetric alternatives with bounded support
 of the radial distributions are well detected by the HZ- and the ${\cal E}$-test. For the case of unbounded support of the radial distribution,
 the strongest test is again HJM. This test is also strongest for the marginally disturbed alternatives MAR$_2$(DIST), where it is just
 outperformed by the PU-test for the disturbance by Exp(1)- and $\chi^2$-random variables. The NM$_d(\vartheta)$-distributions
 are uniformly best detected by HJM, although the power is not very strong, whereas all other tests almost completely fail to detect
 these alternatives. Notably, the S$|\mbox{N}_2|$ alternatives are best detected by $T_{MQ}(\mathbf{f}_1)$.
 Overall, for the chosen alternatives HJM performs best, but it also lacks power especially when the support of the distribution is bounded.
From a robust point of view, the weighted $L^2$ procedures, like DEH$^*$, the HZ-test as well as the energy test ${\cal E}$ perform very well, especially if the focus is on consistency.

In dimensions $d=3$ and $d=5$, one can paint the same picture for the allocation of the best procedures to the alternatives.
Interestingly, the power of the procedures increases compared to the lower-dimensional setting, which appears to be counterintuitive
in view of the curse of dimensionality. Some noticeable phenomena arise: For the $\mathcal{S}^d(\mbox{B}(2,2))$ distribution, some of the tests,
 like HV, HJ and $T_{CS}$, $b_M^{(1)}$, $b_M^{(2)}$ seem to loose power when the sample size is increased.
 An explanation for this behaviour for the latter tests might be that these procedures use an approximation of the maximum
 on the unit sphere, which might be harder to approximate for larger samples. In the case $d=3$, we also observe this behaviour for the HJM-test.
  Interestingly, the HJM-test as well as the PU-test increase the power against NM$_d(\vartheta)$-alternatives in comparison to the case $d=2$,
  whereas the other procedures nearly uniformly fail to distinguish them from the null hypothesis in each dimension considered.

\section{Conclusions and outlook}\label{secconcl}
From a practical point of view, we recommend to use the computationally efficient weighted $L^2$-type procedures.
 like HZ and DEH$^*$,  or the energy test $\mathcal{E}$, since they show a good balance between fast computation time and robust power against many alternatives,
 and they do not exhibit  any particular weakness. If computation time is not an issue we suggest to employ the HJM-test,
  as it outperforms most of the other procedures. Note that by choosing other tuning parameters, the weighted $L^2$-procedures
 are expected to benefit in terms of power against specific alternatives, especially if one is able to choose the tuning parameter in a data dependent way.
  For a first step in this direction for univariate goodness-of-fit tests, see \cite{tenr19}. In general, it would be nice to have explicit solutions
  of the Fredholm integral equation \eqref{int:eq}. For some recent cases in which such integral equations have witnessed explicit
  solutions in the context of goodness-of-fit testing,  see, e.g., Theorem 3.2 of \cite{bata10} or Theorems 3 and 5 of \cite{hr19}. High-dimensional $L^2$-statistics for testing normality have not been considered so far in the literature. The efficient implementation of the tests in the package {\tt mnt} admit first simulations, which indicate that new interesting phenomena arise.

\bibliographystyle{apalike}

\begin{table}[t]
\tiny
\setlength{\tabcolsep}{0.5mm}
\centering
\begin{tabular}{lr|rrrrrrrrrrrrrrrrrr}
 Distribution & $n$ & BHEP & HZ & HV & HJ & HJM & DEH & DEH$^*$ & ${\cal E}$ & $T(\mathbf{f}_1)$ & $T(\mathbf{f}_2)$ & $T_{CS}$ & PU & $b^{(1)}$ & $b^{(2)}$ & $b_{M}^{(1)}$ & $\widetilde{b}^{(2)}$ & $\widetilde{b}^{(1)}$ & $b_{M}^{(2)}$ \\
  \hline
  NMix$(0.5,3,{\rm I}_2)$ & 20  & 18 & 24  & 2 & 3 & 7  & 3  & 16 & 20 & 11 & 19 & 5 & 13 & 2 & 1 & 4 & 2 & 3 & 2 \\
                          & 50  & 64 & 82  & 2 & 2 & 10 & 6  & 67 & 71 & 34 & 51 & 5 & 77 & 2 & 0 & 2 & 0 & 3 & 2 \\
                          & 100 & 99 & 100 & 2 & 2 & 70 & 38 & 99 & 99 & 76 & 88 & 5 & 96 & 2 & 0 & 3 & 0 & 3 & 1 \\
  NMix$(0.79,3,{\rm I}_2)$ & 20  & 42  & 42  & 14 & 10 & 24 & 17 & 34 & 39 & 15 & 13 & 21 & 46 & 18 & 11 & 20 & 11 & 23 & 11 \\
                           & 50  & 94  & 93  & 21 & 6 & 26 & 51 & 89 & 93 & 43 & 23 & 44 & 96 & 52 & 8 & 56 & 8 & 54 & 7 \\
                           & 100 & 100 & 100 & 48 & 5 & 33 & 96 & 100 & 100 & 82 & 50 & 75 & 99 & 91 & 7 & 94 & 6 & 87 & 6 \\
  NMix$(0.9,3,{\rm I}_2)$ & 20 & 38 & 34 & 32 & 26 & 44 & 34 & 37 & 37 & 13 & 22 & 27 & 44 & 34 & 23 & 38 & 24 & 35 & 24 \\
                          & 50 & 83 & 74 & 70 & 38 & 81 & 80 & 83 & 82 & 31 & 62 & 62 & 87 & 87 & 49 & 89 & 51 & 82 & 50 \\
                          & 100 & 99 & 96 & 97 & 45 & 97 & 99 & 99 & 99 & 59 & 95 & 92 & 99 & 100 & 69 & 100 & 72 & 99 & 73 \\
  NMix$(0.5,0,{\rm B}_2)$& 20  & 15 & 14 & 17 & 16 & 34 & 19 & 18 & 16 & 10 & 13 & 12 & 16 & 17 & 18 & 17 & 20 & 16 & 19 \\
                         & 50  & 31 & 31 & 25 & 23 & 60 & 40 & 43 & 34 & 24 & 31 & 15 & 38 & 20 & 32 & 20 & 35 & 18 & 35 \\
                         & 100 & 59 & 61 & 36 & 31 & 84 & 70 & 75 & 63 & 45 & 55 & 16 & 72 & 23 & 52 & 23 & 59 & 20 & 60 \\
  NMix$(0.9,0,{\rm B}_2)$& 20 & 20 & 18 & 27 & 28 & 40 & 28 & 24 & 21 & 10 & 22 & 19 & 22 & 26 & 29 & 26 & 30 & 24 & 30 \\
                         & 50 & 37 & 29 & 54 & 53 & 66 & 55 & 47 & 38 & 15 & 50 & 31 & 43 & 45 & 54 & 45 & 55 & 42 & 57 \\
                         & 100 & 59 & 46 & 78 & 77 & 89 & 78 & 70 & 60 & 24 & 75 & 39 & 66 & 56 & 80 & 56 & 83 & 53 & 85 \\\hline
  t$_1(0,{\rm I}_2)$ & 20 & 96 & 96 & 94 & 94 & 98 & 97 & 97 & 97 & 77 & 97 & 85 & 96 & 92 & 97 & 91 & 97 & 90 & 96 \\
                     & 50 & 100 & 100 & 100 & 100 & 100 & 100 & 100 & 100 & 97 & 100 & 95 & 100 & 99 & 100 & 99 & 100 & 99 & 100 \\
                     & 100 & 100 & 100 & 100 & 100 & 100 & 100 & 100 & 100 & 100 & 100 & 98 & 100 & 100 & 100 & 100 & 100 & 100 & 100 \\
  t$_3(0,{\rm I}_2)$ & 20  & 48 & 45 & 54 & 52 & 67 & 56 & 53 & 49 & 16 & 53 & 36 & 46 & 52 & 59 & 50 & 58 & 48 & 54 \\
                     & 50  & 82 & 78 & 85 & 83 & 95 & 88 & 87 & 83 & 24 & 91 & 55 & 80 & 77 & 92 & 75 & 90 & 72 & 86 \\
                     & 100 & 98 & 97 & 97 & 96 & 100 & 99 & 99 & 98 & 31 & 100 & 96 & 78 & 92 & 99 & 91 & 100 & 87 & 98 \\
  t$_5(0,{\rm I}_2)$ & 20 & 25 & 22 & 32 & 32 & 48 & 34 & 29 & 26 & 9 & 29 & 20 & 24 & 30 & 36 & 31 & 35 & 29 & 31 \\
                     & 50 & 49 & 42 & 58 & 58 & 78 & 63 & 58 & 51 & 11 & 66 & 33 & 46 & 53 & 71 & 52 & 68 & 48 & 62 \\
                     & 100 & 75 & 66 & 81 & 78 & 95 & 83 & 80 & 75 & 12 & 90 & 41 & 68 & 66 & 92 & 64 & 91 & 59 & 86 \\
  t$_{10}(0,{\rm I}_2)$ & 20 & 10 & 9 & 16 & 16 & 25 & 16 & 13 & 11 & 5 & 13 & 10 & 11 & 15 & 17 & 15 & 16 & 15 & 14 \\
                        & 50 & 17 & 13 & 29 & 29 & 48 & 29 & 23 & 18 & 6 & 29 & 15 & 17 & 25 & 32 & 23 & 32 & 22 & 30 \\
                        & 100 & 26 & 20 & 44 & 42 & 69 & 42 & 34 & 27 & 6 & 49 & 18 & 24 & 33 & 57 & 32 & 57 & 30 & 50 \\\hline
  U$^2(0,1)$ & 20  & 12 & 18 & 0 & 0 & 1 & 1 & 10 & 11 & 8 & 33 & 3 & 5 & 0 & 0 & 0 & 0 & 1 & 0 \\
             & 50  & 59 & 68 & 0 & 0 & 7 & 6 & 57 & 52 & 15 & 82 & 2 & 46 & 0 & 0 & 0 & 0 & 0 & 0 \\
             & 100 & 98 & 98 & 0 & 0 & 92 & 79 & 98 & 96 & 27 & 100 & 2 & 96 & 0 & 0 & 0 & 0 & 0 & 0 \\
  LN$^2(0,0.5)$ & 20 & 60 & 55 & 50 & 43 & 55 & 50 & 56 & 57 & 18 & 37 & 52 & 55 & 58 & 40 & 55 & 42 & 58 & 41 \\
                & 50 & 97 & 93 & 92 & 75 & 88 & 93 & 96 & 96 & 45 & 83 & 94 & 95 & 97 & 76 & 95 & 77 & 96 & 73 \\
                & 100 & 100 & 100 & 100 & 93 & 99 & 100 & 100 & 100 & 81 & 99 & 100 & 100 & 100 & 95 & 100 & 96 & 100 & 93 \\
  B$^2(1,2)$ & 20  & 27 & 28 & 5 & 4 & 8 & 6 & 20 & 25 & 10 & 13 & 13 & 18 & 7 & 2 & 7 & 3 & 10 & 3 \\
             & 50  & 81 & 78 & 5 & 1 & 6 & 30 & 73 & 77 & 23 & 41 & 34 & 68 & 21 & 0 & 16 & 1 & 30 & 1 \\
             & 100 & 100 & 99 & 19 & 0 & 17 & 95 & 100 & 99 & 46 & 91 & 77 & 98 & 73 & 0 & 56 & 0 & 74 & 0 \\
  B$^2(2,2)$ & 20  & 5 & 6 & 1 & 1 & 2 & 1 & 3 & 4 & 5 & 12 & 3 & 3 & 1 & 1 & 1 & 1 & 1 & 1 \\
             & 50  & 13 & 17 & 0 & 0 & 0 & 1 & 10 & 10 & 7 & 34 & 2 & 9 & 0 & 0 & 0 & 0 & 0 & 0 \\
             & 100 & 39 & 41 & 0 & 0 & 15 & 7 & 34 & 31 & 9 & 71 & 2 & 27 & 0 & 0 & 0 & 0 & 0 & 0 \\
  P$_{II}^2(0.5,0,1)$ & 20  & 46 & 59 & 0 & 0 & 1 & 1 & 45 & 48 & 20 & 67 & 5 & 24 & 0 & 0 & 1 & 0 & 1 & 0 \\
                      & 50  & 99 & 100 & 0 & 0 & 71 & 69 & 100 & 99 & 52 & 100 & 3 & 98 & 0 & 0 & 0 & 0 & 0 & 0 \\
                      & 100 & 100 & 100 & 0 & 0 & 100 & 100 & 100 & 100 & 85 & 100 & 3 & 100 & 0 & 0 & 0 & 0 & 0 & 0 \\
  P$_{VII}^2(5,0,1)$ & 20 & 20 & 18 & 28 & 27 & 42 & 28 & 25 & 21 & 8 & 23 & 18 & 21 & 28 & 31 & 27 & 30 & 26 & 28 \\
                     & 50 & 38 & 32 & 51 & 49 & 70 & 55 & 48 & 40 & 10 & 55 & 30 & 39 & 43 & 59 & 42 & 59 & 39 & 55 \\
                     & 100 & 63 & 53 & 72 & 71 & 91 & 78 & 72 & 64 & 13 & 82 & 37 & 57 & 56 & 85 & 54 & 85 & 50 & 78 \\
  P$_{VII}^2(10,0,1)$ & 20 & 9 & 8 & 14 & 13 & 25 & 13 & 11 & 10 & 6 & 11 & 9 & 10 & 15 & 17 & 16 & 17 & 14 & 15 \\
                      & 50 & 13 & 11 & 24 & 23 & 41 & 23 & 18 & 14 & 7 & 22 & 13 & 14 & 21 & 28 & 19 & 26 & 20 & 23 \\
                      & 100 & 19 & 14 & 35 & 34 & 57 & 35 & 27 & 20 & 6 & 38 & 16 & 17 & 24 & 44 & 24 & 44 & 22 & 37 \\ \hline
  $\mathcal{S}^2(\mbox{Exp}(1))$ & 20 & 77 & 78 & 68 & 64 & 86 & 75 & 82 & 83 & 32 & 81 & 38 & 67 & 68 & 81 & 64 & 77 & 63 & 70 \\
                                 & 50 & 99 & 100 & 93 & 89 & 99 & 98 & 100 & 100 & 46 & 100 & 49 & 97 & 80 & 99 & 76 & 98 & 72 & 94 \\
                                 & 100 & 100 & 100 & 99 & 98 & 100 & 100 & 100 & 100 & 56 & 100 & 55 & 100 & 91 & 100 & 89 & 100 & 80 & 100 \\
  $\mathcal{S}^2(\mbox{B}(1,2))$ & 20  & 25 & 27 & 14 & 13 & 39 & 22 & 30 & 34 & 11 & 28 & 9 & 18 & 14 & 30 & 12 & 24 & 14 & 17 \\
                                 & 50  & 52 & 66 & 10 & 7 & 60 & 32 & 60 & 70 & 15 & 65 & 6 & 34 & 14 & 44 & 11 & 35 & 12 & 18 \\
                                 & 100 & 84 & 95 & 7 & 2 & 82 & 52 & 90 & 95 & 17 & 93 & 4 & 62 & 9 & 64 & 9 & 57 & 8 & 25 \\
  $\mathcal{S}^2(\mbox{B}(2,2))$ & 20 & 3 & 5 & 0 & 0 & 2 & 1 & 2 & 3 & 4 & 9 & 2 & 3 & 0 & 0 & 1 & 0 & 1 & 0 \\
                                 & 50 & 8 & 10 & 0 & 0 & 0 & 0 & 5 & 6 & 5 & 23 & 1 & 6 & 0 & 0 & 0 & 0 & 0 & 0 \\
                                 & 100 & 21 & 22 & 0 & 0 & 8 & 1 & 14 & 15 & 4 & 54 & 1 & 11 & 0 & 0 & 0 & 0 & 0 & 0 \\
  $\mathcal{S}^2(\chi^2_5)$      & 20 & 16 & 14 & 21 & 21 & 36 & 22 & 19 & 16 & 6 & 18 & 12 & 15 & 19 & 25 & 19 & 24 & 20 & 22 \\
                                 & 50 & 29 & 25 & 38 & 36 & 63 & 42 & 37 & 31 & 7 & 46 & 17 & 26 & 32 & 50 & 29 & 47 & 29 & 40 \\
                                 & 100 & 50 & 43 & 55 & 52 & 85 & 61 & 57 & 52 & 7 & 73 & 20 & 40 & 40 & 76 & 38 & 74 & 35 & 63 \\
  $\mathcal{S}^2(\mbox{LN}(0,0.5))$  & 20 & 12 & 11 & 16 & 16 & 21 & 15 & 14 & 13 & 6 & 13 & 12 & 10 & 16 & 14 & 16 & 15 & 16 & 15 \\
                                     & 50 & 15 & 14 & 30 & 30 & 44 & 29 & 25 & 20 & 7 & 23 & 17 & 15 & 25 & 29 & 25 & 31 & 24 & 30 \\
                                     & 100 & 19 & 22 & 47 & 47 & 64 & 45 & 38 & 31 & 6 & 34 & 22 & 18 & 38 & 50 & 36 & 50 & 34 & 47 \\\hline
  MAR$_2$(Exp(1)) & 20 & 52 & 49 & 41 & 34 & 44 & 39 & 48 & 52 & 16 & 28 & 44 & 57 & 44 & 29 & 48 & 30 & 47 & 29 \\
                  & 50 & 95 & 92 & 81 & 60 & 81 & 86 & 94 & 95 & 40 & 70 & 83 & 97 & 94 & 60 & 96 & 62 & 93 & 62 \\
                  & 100 & 100 & 100 & 99 & 82 & 97 & 100 & 100 & 100 & 75 & 97 & 99 & 100 & 100 & 84 & 100 & 85 & 100 & 86 \\
  MAR$_2(\chi^2_3)$ & 20  & 39 & 35 & 31 & 26 & 36 & 30 & 35 & 37 & 11 & 20 & 32 & 41 & 34 & 23 & 36 & 23 & 36 & 24 \\
                    & 50  & 84 & 77 & 67 & 48 & 67 & 71 & 81 & 83 & 26 & 54 & 70 & 88 & 82 & 48 & 85 & 49 & 80 & 49 \\
                    & 100 & 99 & 98 & 95 & 70 & 89 & 98 & 99 & 100 & 52 & 90 & 95 & 100 & 100 & 68 & 100 & 70 & 98 & 72 \\
  MAR$_2(\chi^2_5)$ & 20 & 25 & 22 & 22 & 20 & 28 & 21 & 22 & 23 & 9 & 15 & 22 & 26 & 25 & 16 & 25 & 16 & 24 & 16 \\
                    & 50 & 61 & 51 & 50 & 35 & 54 & 50 & 57 & 60 & 16 & 38 & 51 & 66 & 64 & 33 & 66 & 35 & 64 & 36 \\
                    & 100 & 93 & 83 & 82 & 52 & 76 & 86 & 91 & 92 & 30 & 71 & 83 & 95 & 93 & 50 & 95 & 51 & 92 & 51 \\
  MAR$_2(t_3)$ & 20 & 23 & 21 & 29 & 28 & 42 & 29 & 27 & 24 & 12 & 24 & 21 & 25 & 27 & 30 & 27 & 31 & 26 & 30 \\
               & 50 & 47 & 41 & 55 & 53 & 73 & 57 & 54 & 48 & 24 & 55 & 33 & 51 & 47 & 58 & 47 & 60 & 44 & 59 \\
                & 100 & 73 & 66 & 78 & 76 & 91 & 81 & 79 & 74 & 42 & 81 & 45 & 76 & 59 & 84 & 58 & 85 & 55 & 85 \\
  MAR$_2(t_5)$ & 20  & 12 & 11 & 16 & 16 & 28 & 17 & 15 & 13 & 7 & 13 & 12 & 13 & 17 & 18 & 17 & 18 & 15 & 19 \\
               & 50  & 21 & 17 & 31 & 30 & 49 & 31 & 27 & 22 & 10 & 29 & 18 & 23 & 28 & 36 & 27 & 36 & 25 & 36 \\
               & 100 & 34 & 27 & 47 & 47 & 70 & 48 & 42 & 35 & 13 & 48 & 24 & 35 & 33 & 52 & 33 & 54 & 31 & 56 \\
  MAR$_2(\Gamma(5,1))$ & 20 & 25 & 22 & 22 & 19 & 30 & 21 & 23 & 24 & 9 & 15 & 22 & 22 & 26 & 18 & 25 & 19 & 28 & 18 \\
                       & 50 & 64 & 53 & 53 & 36 & 55 & 54 & 60 & 62 & 16 & 39 & 60 & 57 & 68 & 36 & 64 & 37 & 68 & 36 \\
                       & 100 & 93 & 83 & 87 & 53 & 79 & 89 & 92 & 92 & 31 & 74 & 93 & 89 & 97 & 52 & 95 & 52 & 96 & 49 \\\hline
  NM$_2(0.2)$ & 20 & 5 & 5 & 6 & 6 & 14 & 6 & 6 & 5 & 5 & 6 & 5 & 5 & 7 & 7 & 7 & 7 & 6 & 7 \\
              & 50 & 5 & 5 & 6 & 6 & 16 & 6 & 6 & 6 & 5 & 6 & 6 & 6 & 7 & 7 & 7 & 7 & 8 & 6 \\
              & 100 & 5 & 5 & 7 & 7 & 17 & 6 & 6 & 6 & 5 & 7 & 6 & 5 & 7 & 7 & 7 & 6 & 6 & 6 \\
  S$|\mbox{N}_2|$ & 20 & 14 & 15 & 11 & 11 & 25 & 15 & 21 & 17 & 33 & 8 & 12 & 13 & 12 & 12 & 12 & 13 & 11 & 12 \\
                  & 50 & 31 & 45 & 16 & 16 & 46 & 50 & 66 & 46 & 87 & 14 & 14 & 33 & 14 & 18 & 14 & 23 & 13 & 24 \\
                  & 100 & 74 & 95 & 23 & 19 & 72 & 96 & 99 & 93 & 100 & 20 & 16 & 66 & 14 & 27 & 15 & 37 & 13 & 38 \\ \hline
  N$_2(\mu_2,\Sigma_{0.5})$     & 20 & 5 & 5 & 5 & 5 & 5 & 5 & 5 & 5 & 5 & 5 & 5 & 5 & 5 & 5 & 5 & 5 & 5 & 5 \\
                                & 50 & 5 & 5 & 5 & 5 & 5 & 5 & 6 & 5 & 5 & 5 & 5 & 5 & 5 & 5 & 5 & 5 & 5 & 5 \\
                                & 100 & 5 & 5 & 5 & 5 & 5 & 5 & 5 & 5 & 5 & 5 & 5 & 5 & 5 & 5 & 5 & 5 & 5 & 5
\end{tabular}
\caption{Empirical rejection rates of the considered tests ($d=2$, $\alpha=0.05$)}\label{tab:d2}
\end{table}

\begin{table}[ht]
\tiny
\setlength{\tabcolsep}{0.5mm}
\centering
\begin{tabular}{lr|rrrrrrrrrrrrrrrrrr}
 Distribution & $n$ & BHEP & HZ & HV & HJ & HJM & DEH & DEH$^*$ & ${\cal E}$ & $T(\mathbf{f}_1)$ & $T(\mathbf{f}_2)$ & $T_{CS}$ & PU & $b^{(1)}$ & $b^{(2)}$ & $b_{M}^{(1)}$ & $\widetilde{b}^{(2)}$ & $\widetilde{b}^{(1)}$ & $b_{M}^{(2)}$ \\
  \hline
NMix$(0.5,3,{\rm I}_3)$    & 20  & 16  & 21  & 3 & 3 & 18 & 3 & 11 & 14 & 13 & 15 & 7 & 34 & 1 & 1 & 3 & 1 & 3 & 2 \\
                           & 50  & 61  & 81  & 3 & 3 & 17 & 5 & 51 & 57 & 39 & 36 & 6 & 77 & 2 & 1 & 2 & 2 & 2 & 3 \\
                           & 100 & 99  & 100 & 3 & 2 & 52 & 30 & 99 & 99 & 87 & 70 & 5 & 86 & 2 & 0 & 2 & 1 & 2 & 2 \\
  NMix$(0.79,3,{\rm I}_3)$ & 20  & 40  & 39 & 13 & 10 & 36 & 16 & 29 & 37 & 16 & 11 & 27 & 75 & 18 & 11 & 18 & 11 & 22 & 10 \\
                           & 50  & 96  & 95 & 14 & 6 & 32 & 46 & 87 & 94 & 49 & 20 & 59 & 96 & 41 & 7 & 47 & 7 & 47 & 6 \\
                           & 100 & 100 & 100 & 25 & 5 & 42 & 95 & 100 & 100 & 91 & 43 & 89 & 98 & 85 & 10 & 93 & 10 & 83 & 6 \\
  NMix$(0.9,3,{\rm I}_3)$  & 20  & 38  & 33 & 33 & 28 & 63 & 35 & 38 & 40 & 14 & 21 & 34    & 72 & 35 & 26 & 40 & 28 & 37 & 28 \\
                           & 50  & 89  & 81 & 66 & 35 & 89 & 84 & 89 & 91 & 40 & 63 & 82    & 96 & 91 & 50 & 94 & 53 & 86 & 52 \\
                           & 100 & 99  & 98 & 96 & 37 & 99 & 100 & 100 & 100 & 79 & 95 & 99 & 99 & 100 & 67 & 100 & 72 & 100 & 73 \\
  NMix$(0.5,0,{\rm B}_3)$  & 20  & 28  & 26 & 28 & 25 & 67 & 36 & 38 & 32 & 17 & 22 & 17 & 54 & 28 & 34 & 25 & 32 & 24 & 24 \\
                           & 50  & 67  & 68 & 45 & 37 & 91 & 77 & 82 & 71 & 45 & 64 & 24 & 83 & 40 & 66 & 35 & 66 & 32 & 50 \\
                           & 100 & 97  & 97 & 63 & 47 & 99 & 98 & 99 & 97 & 82 & 93 & 26 & 98 & 42 & 90 & 39 & 93 & 31 & 79 \\
  NMix$(0.9,0,{\rm B}_3)$  & 20  & 28  & 24 & 43 & 41 & 66 & 42 & 38 & 32 & 11 & 31 & 26 & 55 & 43 & 43 & 43 & 46 & 40 & 44 \\
                           & 50  & 59  & 49 & 80 & 78 & 91 & 80 & 75 & 66 & 20 & 77 & 54 & 80 & 72 & 80 & 73 & 82 & 66 & 81 \\
                           & 100 & 85  & 73 & 96 & 95 & 99 & 96 & 94 & 89 & 34 & 96 & 67 & 94 & 85 & 98 & 84 & 98 & 78 & 97 \\\hline
  t$_1(0,{\rm I}_3)$       & 20  & 99  & 98 & 98 & 97 & 100 & 99 & 99 & 99 & 87 & 99 & 93 & 99 & 98 & 99 & 96 & 99 & 96 & 97 \\
                           & 50  & 100 & 100 & 100 & 100 & 100 & 100 & 100 & 100 & 100 & 100 & 100 & 100 & 100 & 100 & 100 & 100 & 100 & 100 \\
                           & 100 & 100 & 100 & 100 & 100 & 100 & 100 & 100 & 100 & 100 & 100 & 100 & 100 & 100 & 100 & 100 & 100 & 100 & 100 \\
  t$_3(0,{\rm I}_3)$       & 20  & 56 & 52 & 65 & 62 & 88 & 70 & 67 & 62 & 19 & 61 & 43 & 74 & 64 & 71 & 60 & 69 & 56 & 62 \\
                           & 50  & 93 & 91 & 94 & 91 & 99 & 97 & 96 & 95 & 31 & 97 & 74 & 96 & 91 & 98 & 88 & 98 & 84 & 92 \\
                           & 100 & 100 & 100 & 100 & 99 & 100 & 100 & 100 & 100 & 41 & 100 & 87 & 100 & 98 & 100 & 96 & 100 & 94 & 99 \\
  t$_5(0,{\rm I}_3)$       & 20  & 29 & 26 & 42 & 38 & 70 & 44 & 40 & 34 & 10 & 33 & 24 & 52 & 42 & 47 & 37 & 45 & 36 & 38 \\
                           & 50  & 61 & 54 & 74 & 69 & 92 & 80 & 75 & 67 & 13 & 80 & 44 & 77 & 68 & 82 & 64 & 80 & 59 & 72 \\
                           & 100 & 89 & 83 & 92 & 89 & 99 & 96 & 94 & 91 & 14 & 98 & 59 & 92 & 83 & 98 & 79 & 97 & 74 & 92 \\
  t$_{10}(0,{\rm I}_3)$    & 20  & 12 & 10 & 20 & 19 & 46 & 21 & 17 & 14 & 6  & 13 & 11 & 35 & 18 & 21 & 16 & 19 & 17 & 16 \\
                           & 50  & 22 & 17 & 39 & 35 & 68 & 40 & 34 & 27 & 7  & 38 & 19 & 46 & 33 & 46 & 30 & 44 & 30 & 36 \\
                           & 100 & 37 & 28 & 57 & 52 & 86 & 60 & 52 & 43 & 7  & 66 & 25 & 55 & 44 & 73 & 40 & 71 & 35 & 54 \\\hline
  U$^3(0,1)$               & 20  & 11 & 15 & 0 & 0 & 4 & 0 & 4 & 6 & 8 & 37 & 3     & 18 & 0 & 0 & 1 & 0 & 0 & 1 \\
                           & 50  & 58 & 65 & 0 & 0 & 1 & 1 & 32 & 38 & 19 & 88 & 2  & 66 & 0 & 0 & 0 & 0 & 0 & 0 \\
                           & 100 & 98 & 98 & 0 & 0 & 84 & 34 & 95 & 94 & 46 & 100 & 2 & 99 & 0 & 0 & 0 & 0 & 0 & 0 \\
  LN$^3(0,0.5)$            & 20  & 61 & 56 & 55 & 46 & 73 & 54 & 59 & 63 & 20 & 38 & 57 & 78 & 61 & 46 & 53 & 47 & 61 & 43 \\
                           & 50  & 98 & 96 & 93 & 78 & 96 & 96 & 98 & 99 & 54 & 88 & 97 & 99 & 99 & 86 & 96 & 86 & 97 & 79 \\
                           & 100 & 100 & 100 & 100 & 94 & 100 & 100 & 100 & 100 & 90 & 100 & 100 & 100 & 100 & 98 & 100 & 98 & 100 & 96 \\
  B$^3(1,2)$               & 20  & 25 & 26 & 4 & 4 & 15 & 5 & 13 & 20 & 10 & 15 & 13 & 41 & 4 & 2 & 4 & 2 & 8 & 2 \\
                           & 50  & 81 & 78 & 3 & 1 & 11 & 18 & 58 & 75 & 28 & 46 & 33 & 88 & 16 & 0 & 9 & 0 & 30 & 1 \\
                           & 100 & 100 & 99 & 6 & 0 & 20 & 83 & 99 & 100 & 65 & 93 & 75 & 100 & 66 & 0 & 34 & 0 & 72 & 0 \\
  B$^3(2,2)$               & 20  & 5 & 6 & 1 & 1 & 6 & 1 & 2 & 3 & 5 & 17 & 3 & 16 & 1 & 0 & 1 & 1 & 1 & 1 \\
                           & 50  & 14 & 17 & 0 & 0 & 0 & 0 & 4 & 7 & 7 & 41 & 2 & 26 & 0 & 0 & 0 & 0 & 0 & 0 \\
                           & 100 & 40 & 40 & 0 & 0 & 10 & 2 & 20 & 26 & 10 & 79 & 2 & 51 & 0 & 0 & 0 & 0 & 0 & 0 \\
  P$_{II}^3(0.5,0,1)$      & 20  & 38 & 48 & 0 & 0 & 3 & 1 & 20 & 26 & 22 & 69 & 6 & 39 & 0 & 0 & 0 & 0 & 0 & 0 \\
                           & 50  & 99 & 100 & 0 & 0 & 26 & 24 & 97 & 97 & 73 & 100 & 3 & 99 & 0 & 0 & 0 & 0 & 0 & 0 \\
                           & 100 & 100 & 100 & 0 & 0 & 100 & 100 & 100 & 100 & 99 & 100 & 2 & 100 & 0 & 0 & 0 & 0 & 0 & 0 \\
  P$_{VII}^3(5,0,1)$       & 20  & 19 & 17 & 30 & 28 & 58 & 32 & 28 & 23 & 8 & 22 & 18  & 45 & 29 & 32 & 28 & 32 & 27 & 28 \\
                           & 50  & 41 & 34 & 60 & 55 & 84 & 64 & 58 & 47 & 12 & 61 & 35 & 64 & 52 & 70 & 47 & 68 & 46 & 59 \\
                           & 100 & 68 & 58 & 81 & 76 & 97 & 87 & 83 & 73 & 18 & 88 & 48 & 81 & 70 & 92 & 66 & 92 & 59 & 85 \\
  P$_{VII}^3(10,0,1)$      & 20  & 9 & 8 & 14 & 14 & 39 & 14 & 12 & 10 & 6 & 9 & 9 & 31 & 13 & 14 & 13 & 15 & 12 & 13 \\
                           & 50  & 13 & 10 & 26 & 24 & 56 & 27 & 22 & 16 & 6 & 23 & 14 & 36 & 26 & 32 & 22 & 31 & 21 & 26 \\
                           & 100 & 20 & 15 & 40 & 36 & 72 & 41 & 32 & 24 & 7 & 41 & 18 & 41 & 32 & 49 & 30 & 49 & 26 & 40 \\\hline
  $\mathcal{S}^3(\mbox{Exp}(1))$ & 20 & 95 & 95 & 88 & 84 & 99 & 95 & 97 & 97 & 54 & 95 & 63 & 95 & 89 & 97 & 82 & 94 & 81 & 85 \\
                           &  50 & 100 & 100 & 100 & 98 & 100 & 100 & 100 & 100 & 72 & 100 & 80 & 100 & 98 & 100 & 96 & 100 & 92 & 99 \\
                           & 100 & 100 & 100 & 100 & 100 & 100 & 100 & 100 & 100 & 82 & 100 & 86 & 100 & 100 & 100 & 99 & 100 & 95 & 100 \\
  $\mathcal{S}^3(\mbox{B}(1,2))$ & 20 & 63 & 65 & 45 & 37 & 87 & 63 & 72 & 73 & 24 & 64 & 22 & 71 & 50 & 72 & 36 & 60 & 41 & 37 \\
                           & 50  & 97 & 99 & 57 & 40 & 98 & 91 & 98 & 99 & 32 & 99 & 20 & 93 & 55 & 96 & 40 & 91 & 36 & 56 \\
                           & 100 & 100 & 100 & 70 & 35 & 100 & 99 & 100 & 100 & 38 & 100 & 16 & 100 & 58 & 100 & 43 & 100 & 36 & 80 \\
  $\mathcal{S}^3(\mbox{B}(2,2))$  & 20 & 5 & 5 & 3 & 3 & 20 & 4 & 5 & 5 & 5 & 4 & 4 & 21 & 4 & 6 & 3 & 4 & 3 & 4 \\
                           & 50  & 5 & 8 & 0 & 0 & 8 & 3 & 6 & 7 & 5 & 3 & 2 & 21 & 1 & 2 & 1 & 1 & 1 & 0 \\
                           & 100 & 6 & 15 & 0 & 0 & 3 & 3 & 9 & 9 & 5 & 3 & 1 & 19 & 1 & 1 & 0 & 0 & 0 & 0 \\
  $\mathcal{S}^3(\chi^2_5)$  & 20 & 37 & 35 & 44 & 39 & 77 & 50 & 49 & 44 & 10 & 42 & 22 & 59 & 45 & 56 & 37 & 50 & 37 & 39 \\
                           & 50 & 80 & 78 & 73 & 67 & 97 & 86 & 87 & 84 & 11 & 91 & 37 & 83 & 68 & 91 & 58 & 87 & 54 & 71 \\
                           & 100 & 98 & 98 & 91 & 84 & 100 & 98 & 99 & 99 & 11 & 100 & 44 & 97 & 79 & 100 & 70 & 99 & 64 & 93 \\
  $\mathcal{S}^3(\mbox{LN}(0,0.5))$  & 20 & 18 & 15 & 29 & 28 & 56 & 30 & 26 & 22 & 8 & 20 & 17 & 43 & 27 & 30 & 26 & 30 & 26 & 27 \\
                           & 50  & 36 & 28 & 58 & 55 & 85 & 62 & 54 & 43 & 8 & 56 & 32 & 59 & 54 & 69 & 50 & 67 & 47 & 57 \\
                           & 100 & 59 & 45 & 80 & 75 & 96 & 84 & 76 & 66 & 9 & 84 & 41 & 74 & 68 & 91 & 62 & 90 & 57 & 81 \\\hline
  MAR$_3$(Exp(1))          & 20  & 34 & 31 & 30 & 27 & 53 & 28 & 32 & 36 & 12 & 19 & 36 & 69 & 34 & 23 & 37 & 24 & 36 & 26 \\
                           & 50  & 83 & 76 & 66 & 50 & 81 & 70 & 79 & 86 & 30 & 57 & 83 & 97 & 83 & 51 & 87 & 54 & 81 & 55 \\
                           & 100 & 100 & 98 & 94 & 70 & 96 & 97 & 99 & 100 & 64 & 91 & 99 & 100 & 100 & 73 & 100 & 76 & 99 & 79 \\
  MAR$_3(\chi^2_3)$        & 20  & 24 & 21 & 23 & 20 & 47 & 21 & 23 & 26 & 9 & 13 & 25 & 55 & 26 & 19 & 27 & 20 & 28 & 21 \\
                           & 50  & 66 & 56 & 51 & 39 & 68 & 53 & 60 & 69 & 20 & 43 & 67 & 91 & 65 & 35 & 72 & 37 & 64 & 39 \\
                           & 100 & 96 & 90 & 83 & 56 & 89 & 89 & 94 & 97 & 43 & 78 & 95 & 99 & 98 & 60 & 99 & 63 & 96 & 66 \\
  MAR$_3(\chi^2_5)$        & 20  & 16 & 14 & 16 & 15 & 38 & 15 & 15 & 16 & 7 & 10 & 16 & 42 & 18 & 15 & 19 & 15 & 20 & 16 \\
                           & 50  & 43 & 34 & 36 & 27 & 55 & 34 & 38 & 46 & 13 & 28 & 46 & 78 & 48 & 25 & 52 & 26 & 48 & 27 \\
                           & 100 & 79 & 65 & 65 & 40 & 73 & 66 & 73 & 83 & 24 & 56 & 80 & 96 & 82 & 38 & 88 & 40 & 81 & 41 \\
  MAR$_3(t_3)$             & 20  & 17 & 15 & 24 & 24 & 49 & 25 & 23 & 20 & 10 & 17 & 17 & 41 & 22 & 23 & 23 & 25 & 21 & 25 \\
                           & 50  & 35 & 30 & 48 & 47 & 74 & 50 & 45 & 39 & 19 & 45 & 31 & 62 & 39 & 47 & 39 & 49 & 35 & 50 \\
                           & 100 & 61 & 52 & 72 & 69 & 93 & 75 & 71 & 65 & 35 & 72 & 42 & 82 & 58 & 73 & 58 & 77 & 53 & 78 \\
  MAR$_3(t_5)$             & 20  & 9  & 8  & 14 & 13 & 36 & 13 & 12 & 10 & 6  & 9  & 9  & 30 & 11 & 12 & 10 & 12 & 10 & 11 \\
                           & 50  & 15 & 13 & 27 & 25 & 52 & 26 & 23 & 18 & 8  & 22 & 15 & 39 & 23 & 26 & 22 & 26 & 20 & 25 \\
                           & 100 & 24 & 18 & 41 & 40 & 73 & 40 & 34 & 28 & 12 & 38 & 20 & 51 & 32 & 45 & 31 & 48 & 27 & 49 \\
  MAR$_3(\Gamma(5,1))$     & 20  & 25 & 22 & 23 & 20 & 48 & 22 & 23 & 26 & 9 & 14 & 23 & 48 & 27 & 20 & 24 & 21 & 27 & 20 \\
                           & 50  & 66 & 54 & 53 & 37 & 68 & 55 & 60 & 68 & 18 & 42 & 65 & 82 & 70 & 36 & 58 & 38 & 68 & 32 \\
                           & 100 & 96 & 86 & 86 & 54 & 88 & 91 & 94 & 97 & 38 & 80 & 96 & 98 & 98 & 59 & 93 & 59 & 97 & 51 \\\hline
   NM$_3(0.2)$              & 20  & 5 & 5 & 7 & 6 & 24 & 6 & 6 & 5 & 5 & 5 & 6 & 23 & 6 & 6 & 5 & 6 & 5 & 6 \\
                           & 50  & 6 & 6 & 8 & 7 & 27 & 8 & 7 & 6 & 5 & 6 & 7 & 24 & 6 & 6 & 6 & 6 & 5 & 6 \\
                           & 100 & 7 & 6 & 9 & 9 & 36 & 10 & 8 & 7 & 6 & 8 & 7 & 24 & 9 & 12 & 7 & 12 & 8 & 10 \\
  S$|\mbox{N}_3|$          & 20  & 16 & 17 & 14 & 13 & 43 & 20 & 26 & 19 & 36 & 8 & 12 & 39 & 13 & 13 & 11 & 14 & 11 & 12 \\
                           & 50  & 44 & 57 & 21 & 18 & 64 & 63 & 79 & 53 & 93 & 18 & 17 & 68 & 19 & 23 & 18 & 29 & 16 & 23 \\
                           & 100 & 90 & 99 & 30 & 24 & 90 & 99 & 100 & 96 & 100 & 33 & 19 & 84 & 22 & 37 & 21 & 53 & 18 & 40 \\ \hline
  N$_3(\mu_3,\Sigma_{0.5})$     & 20 & 6 & 5 & 5 & 5 & 5 & 5 & 5 & 5 & 5 & 5 & 5 & 5 & 5 & 5 & 5 & 5 & 5 & 5 \\
                                & 50 & 5 & 5 & 5 & 5 & 5 & 5 & 5 & 5 & 5 & 5 & 5 & 5 & 5 & 5 & 5 & 5 & 5 & 5 \\
                                & 100 & 5 & 5 & 5 & 5 & 5 & 5 & 5 & 5 & 5 & 5 & 5 & 5 & 5 & 5 & 5 & 5 & 5 & 5
\end{tabular}
\caption{Empirical rejection rates of the considered tests ($d=3$, $\alpha=0.05$)}\label{tab:d3}
\end{table}

\begin{table}[ht]
\tiny
\setlength{\tabcolsep}{0.5mm}
\centering
\begin{tabular}{lr|rrrrrrrrrrrrrrrrrr}
 Distribution & $n$ & BHEP & HZ & HV & HJ & HJM & DEH & DEH$^*$ & ${\cal E}$ & $T(\mathbf{f}_1)$ & $T(\mathbf{f}_2)$ & $T_{CS}$ & PU & $b^{(1)}$ & $b^{(2)}$ & $b_{M}^{(1)}$ & $\widetilde{b}^{(2)}$ & $\widetilde{b}^{(1)}$ & $b_{M}^{(2)}$ \\
  \hline
NMix$(0.5,3,{\rm I}_5)$         & 20  & 12 & 13 & 4 & 4 & 81 & 3 & 5 & 8 & 10 & 11 & 6   & 40 & 4 & 3 & 3 & 3 & 4 & 3 \\
                                & 50  & 41 & 52 & 3 & 4 & 43 & 3 & 12 & 21 & 27 & 19 & 6 & 45 & 2 & 1 & 2 & 2 & 3 & 3 \\
                                & 100 & 90 & 98 & 3 & 3 & 38 & 9 & 56 & 66 & 70 & 36 & 5 & 46 & 2 & 0 & 2 & 1 & 4 & 2 \\
  NMix$(0.79,3,{\rm I}_5)$      & 20  & 25 & 24 & 12 & 10 & 87 & 12 & 15 & 24 & 13 & 7 & 26 & 46 & 12 & 9 & 14 & 9 & 16 & 9 \\
                                & 50  & 85 & 83 & 9 & 5 & 53 & 22 & 43 & 78 & 36 & 12 & 71 & 56 & 26 & 8 & 32 & 7 & 31 & 5 \\
                                & 100 & 100 & 100 & 11 & 5 & 47 & 66 & 96 & 100 & 82 & 25 & 60 & 95 & 63 & 7 & 86 & 6 & 67 & 3 \\
  NMix$(0.9,3,{\rm I}_5)$       & 20  & 25 & 22 & 32 & 29 & 94 & 27 & 28 & 33 & 13 & 9  & 34 & 48 & 30 & 24 & 38 & 28 & 31 & 28 \\
                                & 50  & 85 & 75 & 50 & 30 & 93 & 65 & 75 & 94 & 41 & 53 & 95 & 59 & 82 & 46 & 95 & 48 & 81 & 49 \\
                                & 100 & 100 & 98 & 77 & 25 & 98 & 98 & 100 & 100 & 86 & 89 & 65 & 89 & 100 & 60 & 100 & 66 & 99 & 68 \\
  NMix$(0.5,0,{\rm B}_5)$       & 20  & 56 & 54 & 53 & 45 & 100 & 72 & 74 & 66 & 30 & 37 & 26  & 68 & 61 & 71 & 36 & 66 & 46 & 39 \\
                                & 50  & 98 & 99 & 78 & 61 & 100 & 99 & 100 & 99 & 79 & 96 & 43 & 87 & 78 & 97 & 55 & 95 & 57 & 64 \\
                                & 100 & 100 & 100 & 93 & 72 & 100 & 100 & 100 & 100 & 99 & 100 & 94 & 30 & 83 & 100 & 66 & 100 & 61 & 90 \\
  NMix$(0.9,0,{\rm B}_5)$       & 20  & 32 & 28 & 62  & 61  & 97  & 56  & 54  & 48 & 15 & 33  & 35 & 60 & 58 & 58 & 57 & 61 & 53 & 58 \\
                                & 50  & 75 & 66 & 96  & 95  & 99  & 95  & 94  & 89 & 26 & 93  & 82 & 84 & 95 & 96 & 92 & 96 & 90 & 95 \\
                                & 100 & 96 & 92 & 100 & 100 & 100 & 100 & 100 & 99 & 45 & 100 & 96 & 95 & 100 & 100 & 99 & 100 & 97 & 100 \\\hline
  t$_1(0,{\rm I}_5)$            & 20  & 99  & 99  & 100 & 99  & 100 & 100 & 100 & 100 & 94  & 99   & 96  & 93 & 100 & 100 & 98 & 100 & 99 & 99 \\
                                & 50  & 100 & 100 & 100 & 100 & 100 & 100 & 100 & 100 & 100 & 100  & 100 & 100 & 100 & 100 & 100 & 100 & 100 & 100 \\
                                & 100 & 100 & 100 & 100 & 100 & 100 & 100 & 100 & 100 & 100 & 100  & 100 & 100 & 100 & 100 & 100 & 100 & 100 & 100 \\
  t$_3(0,{\rm I}_5)$            & 20  & 62  & 59  & 79  & 75  & 100 & 85  & 83  & 76  & 23  & 64   & 50  & 70 & 82  & 86 & 72 & 83 & 73 & 73 \\
                                & 50  & 98  & 97  & 99  & 98  & 100 & 100 & 100 & 99  & 40  & 100  & 89  & 93 & 98  & 100 & 96 & 99 & 95 & 98 \\
                                & 100 & 100 & 100 & 100 & 100 & 100 & 100 & 100 & 100 & 53  & 100  & 98  & 99 & 100 & 100 & 99 & 100 & 99 & 100 \\
  t$_5(0,{\rm I}_5)$            & 20  & 32  & 28  & 54  & 49  & 99  & 59  & 57  & 47  & 12  & 29   & 27  & 58 & 55  & 60 & 49 & 57 & 49 & 50 \\
                                & 50  & 78  & 72  & 89  & 83  & 100 & 95  & 94  & 87  & 15  & 93   & 63  & 80 & 89 & 96 & 81 & 95 & 78 & 85 \\
                                & 100 & 98  & 96  & 99  & 97  & 100 & 100 & 100 & 99  & 19  & 100  & 81  & 92 & 98 & 100 & 94 & 100 & 91 & 98 \\
  t$_{10}(0,{\rm I}_5)$         & 20  & 13  & 12  & 26  & 22  & 95  & 28  & 26  & 20  & 7   & 8    & 12  & 50 & 28 & 31 & 23 & 29 & 23 & 24 \\
                                & 50  & 28 & 23 & 54 & 48 & 93 & 65 & 60 & 44 & 7 & 54 & 28 & 61 & 55 & 67 & 46 & 63 & 44 & 48 \\
                                & 100 & 54 & 43 & 79 & 68 & 98 & 89 & 84 & 69 & 7 & 87 & 39 & 71 & 73 & 92 & 62 & 90 & 55 & 73 \\\hline
  U$^5(0,1)$                    & 20  & 9 & 11 & 0 & 1 & 58 & 0 & 1 & 2 & 7 & 37 & 3 & 33 & 0 & 0 & 1 & 0 & 1 & 1 \\
                                & 50  & 49 & 52 & 0 & 0 & 2 & 0 & 1 & 13 & 20 & 91 & 2 & 44 & 0 & 0 & 0 & 0 & 0 & 0 \\
                                & 100 & 96 & 95 & 0 & 0 & 0 & 0 & 25 & 75 & 57 & 100 & 1 & 69 & 0 & 0 & 0 & 0 & 0 & 0 \\
  LN$^5(0,0.5)$                 & 20  & 54 & 49 & 55 & 48 & 98 & 54 & 56 & 65 & 20 & 26 & 51 & 67 & 62 & 50 & 50 & 50 & 59 & 46 \\
                                & 50  & 98 & 96 & 94 & 81 & 100 & 97 & 98 & 100 & 63 & 91 & 99 & 93 & 100 & 91 & 96 & 91 & 99 & 83 \\
                                & 100 & 100 & 100 & 100 & 96 & 100 & 100 & 100 & 100 & 97 & 100 & 100 & 99 & 100 & 100 & 100 & 100 & 100 & 98 \\
  B$^5(1,2)$                    & 20  & 18 & 19 & 3 & 4 & 78 & 2 & 4 & 12 & 9 & 14 & 11 & 45 & 3 & 2 & 3 & 2 & 6 & 3 \\
                                & 50  & 73 & 68 & 2 & 1 & 26 & 3 & 13 & 62 & 27 & 43 & 24 & 70 & 8 & 0 & 4 & 1 & 23 & 2 \\
                                & 100 & 100 & 98 & 1 & 0 & 12 & 22 & 74 & 99 & 71 & 91 & 65 & 91 & 44 & 0 & 13 & 0 & 65 & 0 \\
  B$^5(2,2)$                    & 20  & 5 & 6 & 1 & 2 & 66 & 1 & 1 & 2 & 5 & 19 & 4 & 35 & 1 & 0 & 1 & 0 & 1 & 1 \\
                                & 50  & 14 & 15 & 0 & 0 & 6 & 0 & 0 & 3 & 8 & 50 & 2 & 37 & 0 & 0 & 0 & 0 & 0 & 0 \\
                                & 100 & 37 & 34 & 0 & 0 & 1 & 0 & 1 & 13 & 13 & 85 & 2 & 46 & 0 & 0 & 0 & 0 & 0 & 0 \\
  P$_{II}^5(0.5,0,1)$           & 20  & 25 & 29 & 0 & 0 & 55 & 0 & 1 & 6 & 17 & 64 & 4 & 35 & 0 & 0 & 0 & 0 & 0 & 0 \\
                                & 50  & 97 & 98 & 0 & 0 & 1 & 0 & 21 & 63 & 71 & 100 & 2 & 63 & 0 & 0 & 0 & 0 & 0 & 0 \\
                                & 100 & 100 & 100 & 0 & 0 & 16 & 24 & 100 & 100 & 100 & 100 & 1 & 86 & 0 & 0 & 0 & 0 & 0 & 0 \\
  P$_{VII}^5(5,0,1)$            & 20  & 16 & 14 & 32 & 29 & 96 & 34 & 32 & 25 & 9 & 12 & 16   & 51 & 32 & 36 & 29 & 35 & 26 & 29 \\
                                & 50  & 39 & 32 & 67 & 61 & 97 & 77 & 73 & 57 & 14 & 66 & 41  & 66 & 64 & 78 & 56 & 76 & 54 & 63 \\
                                & 100 & 71 & 60 & 88 & 82 & 100 & 95 & 94 & 83 & 24 & 94 & 58 & 79 & 84 & 96 & 75 & 95 & 70 & 87 \\
  P$_{VII}^5(10,0,1)$           & 20  & 8 & 7 & 14 & 13 & 92 & 15 & 14 & 12 & 6 & 4 & 8     & 44 & 15 & 16 & 12 & 16 & 12 & 13 \\
                                & 50  & 11 & 9 & 28 & 25 & 80 & 32 & 29 & 19 & 6 & 21 & 15  & 51 & 27 & 37 & 22 & 34 & 20 & 26 \\
                                & 100 & 18 & 14 & 45 & 38 & 87 & 54 & 46 & 28 & 8 & 45 & 19 & 54 & 35 & 56 & 31 & 55 & 28 & 40 \\\hline
  $\mathcal{S}^5(\mbox{Exp}(1))$ & 20 & 99 & 99 & 99 & 97 & 100 & 100 & 100 & 100 & 80 & 99 & 83 & 91 & 99 & 100 & 95 & 100 & 95 & 94 \\
                                & 50  & 100 & 100 & 100 & 100 & 100 & 100 & 100 & 100 & 93 & 100 & 98 & 99 & 100 & 100 & 100 & 100 & 99 & 100 \\
                                & 100 & 100 & 100 & 100 & 100 & 100 & 100 & 100 & 100 & 97 & 100 & 100 & 100 & 100 & 100 & 100 & 100 & 100 & 100 \\
  $\mathcal{S}^5(\mbox{B}(1,2))$ & 20 & 94 & 94 & 86 & 75 & 100 & 97 & 98 & 97 & 52 & 90 & 47      & 82 & 93 & 98 & 70 & 95 & 77 & 72 \\
                                & 50  & 100 & 100 & 98 & 90 & 100 & 100 & 100 & 100 & 62 & 100 & 63 & 96 & 98 & 100 & 83 & 100 & 81 & 92 \\
                                & 100 & 100 & 100 & 100 & 96 & 100 & 100 & 100 & 100 & 69 & 100 & 64 & 99 & 99 & 100 & 92 & 100 & 85 & 100 \\
  $\mathcal{S}^5(\mbox{B}(2,2))$ & 20 & 26 & 27 & 22 & 16 & 96 & 42 & 43 & 34 & 9 & 16 & 9 & 53 & 29 & 43 & 11 & 33 & 18 & 14 \\
                                & 50  & 65 & 70 & 14 & 6 & 90 & 62 & 73 & 67 & 8 & 68 & 9 & 61 & 26 & 66 & 9 & 40 & 13 & 9 \\
                                & 100 & 94 & 97 & 6 & 1 & 95 & 82 & 95 & 93 & 7 & 95 & 8 & 68 & 22 & 88 & 8 & 69 & 11 & 5 \\
  $\mathcal{S}^5(\chi^2_5)$     & 20  & 72 & 70 & 77 & 67 & 100 & 89 & 88 & 82 & 22 & 69 & 39       & 72 & 78 & 88 & 61 & 82 & 64 & 62 \\
                                & 50  & 100 & 100 & 98 & 93 & 100 & 100 & 100 & 100 & 23 & 100 & 72 & 91 & 98 & 100 & 89 & 100 & 87 & 94 \\
                                & 100 & 100 & 100 & 100 & 99 & 100 & 100 & 100 & 100 & 24 & 100 & 84 & 98 & 100 & 100 & 98 & 100 & 93 & 100 \\
  $\mathcal{S}^5(\mbox{LN}(0,0.5))$ & 20 & 33 & 30 & 55 & 49 & 99 & 62 & 60 & 49 & 11 & 31 & 26   & 62 & 57 & 66 & 45 & 61 & 47 & 48 \\
                                & 50  & 82 & 77 & 91 & 83 & 100 & 97 & 96 & 90 & 14 & 95 & 61     & 81 & 88 & 97 & 80 & 94 & 76 & 83 \\
                                & 100 & 99 & 98 & 99 & 97 & 100 & 100 & 100 & 100 & 15 & 100 & 78 & 93 & 98 & 100 & 93 & 100 & 88 & 98 \\\hline
  MAR$_5$(Exp(1))               & 20 & 16 & 14 & 19 & 19 & 90 & 16 & 17 & 20 & 9 & 6 & 22    & 50 & 19 & 16 & 22 & 17 & 22 & 18 \\
                                & 50 & 47 & 39 & 45 & 37 & 83 & 42 & 44 & 62 & 17 & 36 & 72  & 67 & 59 & 35 & 67 & 39 & 59 & 39 \\
                                & 100 & 88 & 76 & 76 & 55 & 93 & 76 & 80 & 96 & 40 & 72 & 98 & 74 & 94 & 56 & 97 & 64 & 91 & 66 \\
  MAR$_5(\chi^2_3)$             & 20 & 11 & 11 & 15 & 14 & 88 & 12 & 13 & 14 & 7 & 5 & 15    & 48 & 14 & 11 & 17 & 13 & 15 & 14 \\
                                & 50 & 32 & 26 & 34 & 28 & 80 & 30 & 31 & 45 & 12 & 25 & 52  & 64 & 45 & 28 & 56 & 30 & 45 & 34 \\
                                & 100 & 71 & 56 & 59 & 41 & 88 & 58 & 61 & 84 & 25 & 54 & 91 & 71 & 84 & 40 & 94 & 47 & 80 & 53 \\
  MAR$_5(\chi^2_5)$             & 20 & 9 & 8 & 11 & 11 & 88 & 9 & 9 & 10 & 6 & 5 & 10        & 44 & 11 & 10 & 13 & 11 & 12 & 11 \\
                                & 50 & 21 & 16 & 24 & 19 & 68 & 19 & 19 & 28 & 9 & 16 & 31   & 58 & 28 & 15 & 35 & 17 & 30 & 21 \\
                                & 100 & 46 & 33 & 41 & 29 & 75 & 37 & 38 & 61 & 16 & 35 & 71 & 68 & 60 & 27 & 78 & 33 & 61 & 37 \\
  MAR$_5(t_3)$                  & 20 & 9 & 9 & 18 & 17 & 90 & 15 & 15 & 14 & 8 & 7 & 12      & 45 & 17 & 15 & 21 & 18 & 16 & 21 \\
                                & 50 & 20 & 16 & 39 & 38 & 81 & 38 & 37 & 30 & 13 & 31 & 25  & 54 & 31 & 34 & 36 & 38 & 29 & 40 \\
                                & 100 & 38 & 30 & 62 & 60 & 91 & 63 & 59 & 51 & 22 & 58 & 38 & 61 & 50 & 62 & 56 & 67 & 47 & 70 \\
  MAR$_5(t_5)$                  & 20 & 6 & 6 & 11 & 10 & 86 & 9 & 8 & 7 & 5 & 5 & 7          & 42 & 9 & 8 & 10 & 9 & 8 & 9 \\
                                & 50 & 9 & 8 & 20 & 19 & 69 & 19 & 17 & 13 & 6 & 13 & 12     & 46 & 18 & 20 & 18 & 22 & 16 & 22 \\
                                & 100 & 13 & 10 & 32 & 30 & 73 & 31 & 27 & 20  & 8 & 24 & 17 & 51 & 24 & 32 & 28 & 35 & 21 & 37 \\
  MAR$_5(\Gamma(5,1))$          & 20 & 19 & 17 & 21 & 18 & 93 & 18 & 19 & 23 & 9 & 7 & 18    & 52 & 20 & 17 & 18 & 18 & 22 & 17 \\
                                & 50 & 59 & 47 & 50 & 36 & 87 & 50 & 53 & 73 & 19 & 39 & 62  & 79 & 73 & 42 & 54 & 42 & 71 & 37 \\
                                & 100 & 95 & 84 & 83 & 54 & 95 & 88 & 91 & 99 & 46 & 80 & 97 & 91 & 98 & 67 & 90 & 68 & 97 & 53 \\\hline
  NM$_5(0.2)$                   & 20 & 6 & 6 & 9 & 9 & 87 & 10 & 9 & 8 & 6 & 4 & 7           & 42  & 8  & 9  & 8 & 9 & 7 & 7 \\
                                & 50 & 10 & 9 & 12 & 11 & 68 & 16 & 16 & 12 & 10 & 8 & 9     & 45  & 12 & 15 & 12 & 14 & 11 & 13 \\
                                & 100 & 13 & 12 & 17 & 14 & 75 & 26 & 25 & 16 & 16 & 15 & 10 & 47  & 14 & 22 & 14 & 26 & 11 & 22 \\
  S$|\mbox{N}_5|$               & 20 & 16 & 16 & 15 & 13 & 95 & 20 & 23 & 18 & 26 & 4 & 10   & 48 & 16 & 17 & 12 & 16 & 12 & 11 \\
                                & 50 & 41 & 47 & 24 & 19 & 88 & 54 & 67 & 42 & 86 & 19 & 17  & 57 & 26 & 33 & 17 & 38 & 18 & 20 \\
                                & 100 & 87 & 95 & 36 & 26 & 95 & 97 & 99 & 82 & 100 & 40 & 21 & 65 & 28 & 53 & 22 & 68 & 19 & 35 \\ \hline
  N$_5(\mu_5,\Sigma_{0.5})$     & 20 & 5 & 5 & 5 & 5 & 5 & 5 & 5 & 5 & 5 & 5 & 5 & 5 & 5 & 5 & 5 & 5 & 5 & 5 \\
                                & 50 & 5 & 5 & 5 & 5 & 5 & 5 & 5 & 5 & 5 & 5 & 5 & 5 & 5 & 5 & 5 & 5 & 5 & 5 \\
                                & 100 & 5 & 5 & 5 & 4 & 5 & 6 & 5 & 5 & 5 & 5 & 5 & 5 & 5 & 5 & 5 & 5 & 5 & 5
\end{tabular}
\caption{Empirical rejection rates of the considered tests ($d=5$, $\alpha=0.05$)}\label{tab:d5}
\end{table}

\end{document}